\documentclass[a4paper]{amsart}

\pdfoutput=1

\usepackage[a4paper,inner=3cm,outer=3cm,top=3.3cm,bottom=3cm]{geometry}
\usepackage[colorlinks=true,citecolor=blue, urlcolor=magenta, breaklinks=true]{hyperref}

\usepackage[english]{babel}
\usepackage{geometry}
\usepackage{standalone}
\usepackage{amsfonts,amsmath,amssymb,amsthm,mathrsfs}
\usepackage{caption}
\usepackage[labelfont=rm]{subcaption}
\usepackage{xspace}
\usepackage{microtype}
\usepackage{enumerate}
\usepackage{booktabs}
\usepackage{float}
\usepackage{rotating}
\usepackage{minipage}
\usepackage{multirow}
\usepackage{array}
\usepackage{mathtools}
\usepackage{lscape}
\usepackage{tikz}
\usetikzlibrary{arrows}
\usepackage{array,boldline, makecell, booktabs}
\usepackage{pgf}
\usetikzlibrary{decorations.markings}

\microtypesetup{tracking,kerning,spacing}
\microtypecontext{spacing=nonfrench}

\newcommand\polymake{{\tt polymake}\xspace}

\DeclareSymbolFontAlphabet{\mathbb}{AMSb}

\newcommand\arXiv[1]{\href{http://arxiv.org/abs/#1}{\texttt{arXiv:#1}}}

\newcommand{\pil}{\rightarrow}
\newcommand{\nn}{\mathbb{N}}
\newcommand{\zz}{\mathbb{Z}}
\newcommand{\rr}{\mathbb{R}}
\newcommand{\La}{\Lambda}

\newcommand{\mc}{\text{mc}}
\newcommand{\cc}{\mathcal{C}}
\newcommand{\ks}{\mathcal{S}}
\newcommand{\ora}{\overrightarrow}
\newcommand{\kr}{\mathcal{R}}
\newcommand{\xx}{\mathfrak{X}}
\newcommand{\LL}{\mathfrak{L}}
\newcommand{\sub}{\subseteq}

\newcommand{\smallskeleton}{
\coordinate (A) at (0,0);
\draw [thick] (-0.5,0.5)--(A);
\draw [thick] (-0.5,-0.5)--(A);
\draw [thick] (A) --(1,0);
\draw [thick] (1,0) --(1.5,0.5);
\draw [thick] (1,0) --(1.5,-0.5);
}
\newcommand{\decskeleton}{
\draw [thick] (-0.5,0.5)--(0,0);
\draw [thick] (-0.5,-0.5)--(0,0);
\draw [thick] (0,0) --(1,0);
\draw [thick] (1,0) --(1.5,0.5);
\draw [thick] (1,0) --(1.5,-0.5);
 \node [above left] at (-0.5,0.5) {$j$};
 \node [below left] at (-0.5,-0.5) {$i$};
 \node [above right] at (1.5,0.5) {$l$};
 \node [below right] at (1.5,-0.5) {$k$};
 }
 \newcommand{\dmA}{
 \draw [ultra thick] (-1,0)--(0,0);
\draw [ultra thick] (-1,0)--(-1,1);
\draw [ultra thick] (-1,1) --(0,1);
\draw [thick] (0,1) --(0,0);
\draw [thick] (-1,0) --(0,1);
\draw [thick] (0,0) --(1,-0.3);
\draw [ultra thick] (1.3,0.3) --(1,-0.3);
\draw [thick] (1.3,0.3) --(0,1);
\draw [thick] (0,1) --(1,-0.3);
\draw [thick, dashed] (0,0) --(1.3,0.3);
 \node [below left] at (-1,0) {$A$};
\node [above left] at (-1,1) {$B$};
\node [below] at (0,0) {$C$};
\node [ above right] at (0,1) {$D$};
\node [right] at (1,-0.3) {$E$};
\node [ right] at (1.3,0.3) {$F$};
}
\newcommand{\dmB}{
\draw [ultra thick] (-1,0)--(0,0);
\draw [thick] (0,1) --(0,0);
\draw [ultra thick] (-1,0) --(0,1);
\draw [thick] (0,0) --(1,-0.3);
\draw [thick] (1.3,0.3) --(1,-0.3);
\draw [thick] (1.3,0.3) --(0,1);
\draw [thick] (0,1) --(1,-0.3);
\draw [ultra thick] (1,-0.3) --(2.3,0);
\draw [ultra thick] (1.3,0.3) --(2.3,0);
\draw [thick, dashed] (0,0) --(1.3,0.3);
 \node [below left] at (-1,0) {$A$};
\node [below] at (0,0) {$B$};
\node [ above right] at (0,1) {$C$};
\node [right] at (2.3,0) {$F$};
\node [below right] at (1,-0.3) {$D$};
\node [ above right] at (1.3,0.3) {$E$};
}
\newcommand{\dmC}{
\draw [ultra thick] (-1,0)--(0,0);
\draw [thick] (0,1) --(0,0);
\draw [ultra thick] (-1,0) --(0,1);
\draw [thick] (0,0) --(1,-0.3);
\draw [ultra thick] (1.3,0.3) --(1,-0.3);
\draw [thick] (1.3,0.3) --(0,1);
\draw [thick] (0,1) --(1,-0.3);
\draw [thick] (1,-0.3) --(2.3,0);
\draw [thick, dashed] (1.3,0.3) --(2.3,0);
\draw [thick, dashed] (0,0) --(1.3,0.3);
\draw [ultra thick](2.3,0)--(2.3,1);
\draw [thick](1.3,0.3)--(2.3,1);
\draw [thick](1,-0.3)--(2.3,1);
 \node [below left] at (-1,0) {$A$};
\node [below] at (0,0) {$B$};
\node [ above right] at (0,1) {$C$};
\node [right] at (2.3,0) {$F$};
\node [below right] at (1,-0.3) {$D$};
\node [ above] at (1.3,0.3) {$E$};
\node [ above right] at (2.3,1) {$G$};
}
\newcommand{\dmD}{
\draw [thick] (0,0) --(0.5,1);
\draw [ultra thick] (0.5,1) --(1,0);
\draw [ultra thick] (0,0) --(1,0);
\draw[thick](0,0)--(-1,-0.3);
\draw[thick](1,0)--(2,-0.3);
\draw[ultra thick](-1.3,0.3)--(-1,-0.3);
\draw[ultra thick](2.3,0.3)--(2,-0.3);
\draw[thick](-1.3,0.3)--(0.5,1);
\draw[thick](2.3,0.3)--(0.5,1);
\draw[thick](0.5,1)--(-1,-0.3);
\draw[thick](0.5,1)--(2,-0.3);
\draw[thick, dashed](-1.3,0.3)--(0,0);
\draw[thick,dashed](2.3,0.3)--(1,0);
 \node [below left] at (-1,-0.3) {$A$};
  \node [above left] at (-1.3,0.3) {$B$};
\node [below] at (0,0) {$C$};
\node [ above] at (0.5,1) {$D$};
\node [below] at (1,0) {$E$};
\node [below right] at (2,-0.3) {$F$};
\node [ above right] at (2.3,0.3) {$G$};
}
\newcommand{\dmE}{
\draw [thick] (-0.5,1) --(-0.5,0);
\draw [thick] (-0.5,0) --(0.5,-0.3);
\draw [ultra thick] (1,-0.3) --(0.5,-0.3);
\draw [thick] (1,-0.3) --(-0.5,1);
\draw [thick] (-0.5,1) --(0.5,-0.3);
\draw [thick] (0.5,1.5) --(-0.5,1);
\draw [thick] (1,1.5) --(-0.5,1);
\draw [thick,dashed] (0.5,1.5) --(-0.5,0);
\draw [thick] (1,1.5) --(-0.5,0);
\draw [ultra thick] (1,1.5) --(0.5,1.5);
\draw [ultra thick] (-1.5,0.5) --(-0.5,0);
\draw [ultra thick] (-1.5,0.5) --(-0.5,1);
\draw [thick, dashed] (-0.5,0) --(1,-0.3);
\node [left] at (-1.5,0.5) {$A$};
\node [below] at (-0.5,0) {$B$};
\node [above] at (-0.5,1) {$C$};
\node [below] at (0.5,-0.3) {$D$};
\node [below] at (1,-0.3) {$E$};
\node [above] at (0.5,1.5) {$F$};
\node [above right] at (1,1.5) {$G$};
}
\newcommand{\dmF}{
\draw [thick] (-1.5,-0.3) --(-0.5,0);
\draw [ultra thick](-1.5,-0.3)--(-1.7,0.3);
\draw [thick](-1.5,-0.3)--(-0.5,1);
\draw [thick](-1.7,0.3)--(-0.5,1);
\draw [ultra thick] (-0.5,1) --(-0.5,0);
\draw [thick,dashed] (-1.7,0.3) --(-0.5,0);
\draw [thick] (-0.5,0) --(0.5,-0.3);
\draw [ultra thick] (0.8,0.3) --(0.5,-0.3);
\draw [thick] (0.8,0.3) --(-0.5,1);
\draw [thick] (-0.5,1) --(0.5,-0.3);
\draw [thick] (0.5,-0.3) --(1.8,0);
\draw [thick, dashed] (0.8,0.3) --(1.8,0);
\draw [thick, dashed] (-0.5,0) --(0.8,0.3);
\draw [ultra thick] (1.8,1) --(1.8,0);
\draw [thick] (0.5,-0.3) --(1.8,1);
\draw [thick] (0.8,0.3) --(1.8,1);
\node [left] at (-1.5,-0.3) {$A$};
\node [left] at (-1.7,0.3) {$B$};
\node [below] at (-0.5,0) {$C$};
\node [above] at (-0.5,1) {$D$};
\node [below] at (0.5,-0.3) {$E$};
\node [above] at (0.8,0.3) {$F$};
\node [ right] at (1.8,0) {$G$};
\node [right] at (1.8,1) {$H$};
}
\newcommand{\dmG}{
\draw [thick] (-0.5,1) --(-0.5,0);
\draw [thick] (-0.5,0) --(0.5,-0.3);
\draw [ultra thick] (0.8,0.3) --(0.5,-0.3);
\draw [thick] (0.8,0.3) --(-0.5,1);
\draw [thick] (-0.5,1) --(0.5,-0.3);
\draw [thick] (0.5,-0.3) --(1.8,0);
\draw [thick, dashed] (0.8,0.3) --(1.8,0);
\draw [thick, dashed] (-0.5,0) --(0.8,0.3);
\draw [ultra thick] (1.8,1) --(1.8,0);
\draw [thick] (0.5,-0.3) --(1.8,1);
\draw [thick] (0.8,0.3) --(1.8,1);
\node [left] at (-0.5,0) {$A$};
\node [ left] at (-0.5,1) {$B$};
\node [below] at (0.5,-0.3) {$C$};
\node [above] at (0.8,0.3) {$D$};
\node [ right] at (1.8,0) {$E$};
\node [right] at (1.8,1) {$F$};
}
\newcommand{\dmH}{
\draw [thick] (-0.5,1) --(-0.5,0);
\draw [thick] (-0.5,0) --(0.5,-0.3);
\draw [thick] (0.8,0.3) --(0.5,-0.3);
\draw [thick] (0.8,0.3) --(-0.5,1);
\draw [thick] (-0.5,1) --(0.5,-0.3);
\draw [ultra thick] (0.5,-0.3) --(1.8,0);
\draw [ultra thick] (0.8,0.3) --(1.8,0);
\draw [thick, dashed] (-0.5,0) --(0.8,0.3);
\node [left] at (-0.5,0) {$A$};
\node [ left] at (-0.5,1) {$B$};
\node [below] at (0.5,-0.3) {$C$};
\node [above right] at (0.8,0.3) {$D$};
\node [ right] at (1.8,0) {$E$};
}
\newcommand{\dmI}{
\draw [thick] (-0.5,1) --(-0.5,0);
\draw [thick] (-0.5,0) --(0.5,-0.3);
\draw [ultra thick] (0.8,0.3) --(0.5,-0.3);
\draw [thick] (0.8,0.3) --(-0.5,1);
\draw [thick] (-0.5,1) --(0.5,-0.3);
\draw [thick, dashed] (-0.5,0) --(0.8,0.3);
\node [left] at (-0.5,0) {$A$};
\node [ left] at (-0.5,1) {$B$};
\node [below] at (0.5,-0.3) {$C$};
\node [above right] at (0.8,0.3) {$D$};
}
\newcommand{\dmJ}{
 \draw [ultra thick] (0,0)--(1,0);
\draw [thick,dashed] (-1,1)--(0,0);
\draw [thick](-1,1)--(1,0);
\draw[thick](-1.5,-0.5)--(1,0);
\draw [thick,dashed] (-1.5,-0.5) --(0,0);
\draw [thick,dashed] (-2,0.5) --(0,0);
\draw [thick,dashed] (-2,0.5) --(1,0);
\draw [thick] (-1.5,-0.5)--(-1,1);
\draw [thick](-2,0.5)--(-1,1);
\draw[ultra thick](-1.5,-0.5)--(-2,0.5);
 \node  [above] at (-1,1) {$A$};
 \node [right] at (1,0) {$B$};
 \node [below] at (0,0) {$C$};
 \node [left] at (-2,0.5) {$D$};
 \node [left] at (-1.5,-0.5) {$E$};
}
\DeclareMathOperator{\inter}{int}

\DeclareMathOperator{\rank}{rank}

\theoremstyle{plain}
    \newtheorem{theorem}{Theorem}

    \newtheorem{lemma}[theorem]{Lemma}
    \newtheorem{proposition}[theorem]{Proposition}
      
      \newtheorem{question}[theorem]{Question}
\theoremstyle{definition}
    \newtheorem{remark}[theorem]{Remark}
    
    \newtheorem{definition}[theorem]{Definition}

\author[M.\,Panizzut]{Marta Panizzut}
\address[M.\,Panizzut]{TU Berlin\\Germany}
\email{panizzut@math.tu-berlin.de}

\author[M. D. Vigeland]{Magnus Dehli Vigeland}
\address[M. D. Vigeland]{Oslo, Norway}
\email{m.d.vigeland@medisin.uio.no}

\keywords{Tropical cubic surfaces, tropical lines}

\title{Tropical Lines on Cubic Surfaces}

\begin{document}

\begin{abstract}
Given a tropical line $L$ and a smooth tropical surface $X$, we look at the position of $L$ on $X$. We introduce its  primal and dual motif which are respectively a decorated graph and a subcomplex of the dual triangulation of $X$. They encode the combinatorial  position of $L$ on $X$.  We classify all possible motifs of tropical lines on general smooth tropical surfaces. This classification allows  to give an upper bound for the number of tropical lines on a general smooth tropical surface with a given subdivision. We  focus  in particular on surfaces of degree three.  As a concrete example, we look at tropical cubic surfaces dual to a fixed honeycomb triangulation, showing that a general surface contains exactly $27$ tropical lines.
\end{abstract}

\maketitle
\section{Introduction}
Arthur Cayley stated in his 1849 correspondence with George Salmon  that  a general cubic surface in three-space contains a finite number of lines. In subsequent manuscripts they proved one of the most celebrated result in nineteenth century algebraic geometry: any smooth cubic surface contains exactly $27$ distinct lines. Few years later, in 1854, Ludwig Schl\"afli determined how the lines intersect. A \emph{double-six} is a pair of sets  of $6$ skew lines. He  showed that a general cubic surface has $36$ double-sixes. 

Since the early days of tropical geometry \cite{MaclaganSturmfels15}, it was  a recurring question whether the Cayley-Salmon Theorem has a tropical analogue. The statement that any smooth tropical cubic surface should contain exactly $27$~lines turned out to be false, as there are tropical surfaces containing infinite families of  tropical lines \cite{Vigeland10}. Therefore, it was far from obvious what the correct formulation of the tropical analogue is. 

The story of this paper begins about twelve years ago, when the second author developed a systematic approach to classify tropical lines on smooth tropical surfaces. This led to the discovery of smooth tropical surfaces of any degree with infinitely many tropical lines \cite{Vigeland10}. Furthermore, a classification of tropical lines on smooth tropical cubic surfaces remained in an unpublished manuscript \cite{Vigeland08}.  In 2018, the first author revisited the material in \cite{Vigeland08}. The outcome is this joint publication presented here.
 
 Tropical lines and tropical surfaces are unbounded rational polyhedral complexes of dimension one and two respectively. 
 We present a way of displaying information of the combinatorial position of lines by decorating edges and vertices of their underlying graphs. A key ingredient is the duality between combinatorial types of smooth tropical surfaces and regular unimodular triangulations of their Newton polytopes.  The main idea of \cite{Vigeland08} was to classify all possible combinatorial positions that a tropical line can have on a smooth tropical surface, and to analyze which of these could occur on cubic surfaces.

Hampe and Joswing \cite[Section 6.2]{HampeJoswig17} found  a gap in the classification presented in  \cite{Vigeland08}.  They exhibited an example of a general smooth tropical cubic surface containing $26$ isolated lines and $3$ families. One of these families has  a  combinatorial position not reported in \cite{Vigeland08}. The  present paper fully completes the classification.  Borrowing a term from computational biology, we introduce \emph{motifs} of lines on tropical surfaces which encode combinatorial positions. 
 We establish the following result for tropical cubic surfaces:  
\begin{theorem}\label{theo:classif}
Let $X$ be a general smooth tropical cubic surface. Non-degenerate tropical lines on $X$ have the following 10 motifs illustrated  in Table \ref{tab:gen_comblines}: 
\begin{itemize}
\item \emph{3A}, \emph{3B}, \emph{3C}, \emph{3D}, \emph{3E}, \emph{3F}, \emph{ 3G} and \emph{3H} for  isolated lines,   
\item \emph{3I} and \emph{3J} for families of lines. 
\end{itemize}
\end{theorem}
In the last section of this article, we exhaustively look at one example. We  study  smooth tropical cubic surfaces (not necessarily general)  dual to the honeycomb triangulation in Figure \ref{subex}. This is an interesting triangulation to consider due to its nice symmetries. Motivated by this example, the second author posed two conjectures in the original manuscript: Firstly, that any general smooth cubic surface contains exactly $27$~tropical lines; and secondly, that for any smooth cubic surface the number of families plus the number of isolated lines equal~27. The example provided by Hampe and Joswig gives a counterexample to both conjectures, as the surface contains $26$~isolated lines and $3$~families. They were able to analyze the surface and their lines thanks to the software \texttt{a-tint} \cite{atint}, which is  an extension of \texttt{polymake}  \cite{polymake} for tropical intersection theory. 

The polyhedral structure of tropical varieties makes tropical geometry well suited for computations. There are numerous software packages to deal with tropical geometry, including \texttt{gfan} \cite{gfan}, \polymake \cite{polymake}, \texttt{Singular} \cite{DGPS}, and \texttt{TropLi} \cite{tropli}.  This makes tropical geometry and computational geometry inextricably connected. We strongly believe that this connection holds the  key for further investigations of the geometry of lines on tropical cubic surfaces. Thanks to their \emph{ down-flip reverse search} algorithm  for parallel enumeration of triangulations implemented in \texttt{MPTOPCOM} \cite{mptopcom}, Jordan, Joswig and Kastner computed the  14\,373\,645 combinatorial types of smooth tropical cubic surfaces. This result brings up exciting computational challenges in the pursuit of providing a full census of tropical lines on smooth tropical cubic surfaces.  Through the analysis carried out in this article, we are able to determine necessary conditions for   motifs in a regular unimodular triangulation, see Table \ref{tab:gen_comblines}. However, they  are not sufficient: a motif may or may not correspond to a line on a surface dual to the triangulation. As the study of the honeycomb triangulation  highlights, there are more conditions  that needs to be imposed to the coefficients of the surface to  assure realizability.  These determine further subdivisions of the secondary cone.  We believe that carrying out the same analysis for all unimodular triangulations  is the next step in deepening the classification. 

Smooth cubic algebraic surfaces are Del Pezzo surfaces of degree three.  In  \cite{RSS16},  Ren, Shaw and Sturmfels studied the intrinsic tropicalization of very affine varieties over valued field obtained from Del Pezzo surfaces by removing the  $27$ lines. Each tropical line has ten markings given by the intersection with ten other lines. These markings determine different combinatorial types that they classified in Table 1. This approach is different from the one we are considering in this paper. As the authors say in the introduction, \lq\lq there are different tropical models of a single classical variety, and the choice of model depends on what structure one wants revealed''.  However, an interesting question is to create a tropical connection between these two models and to understand what happens to the incidence structure of the lines on the tropical Del Pezzo surface when we  consider an embedding in $\mathbb{R}^3$. 

Finally, we want to mention the lifting problem of a pair given by a tropical surface~$X$ and a line~$L$ on it. This problem deals with understanding whether there exist an algebraic surface and a line on it which tropicalize to $X$ and $L$, respectively. We recently learned about the work of Geiger \cite{Geiger19} who in her Master Thesis also looks at the classification of two-points families on general smooth cubic surfaces. She  presents a first attempt to solve the lifting problem for infinite families of lines. In particular, she proves that the problem has a negative answer for non-degenerates lines in families of motif 3I over fields with residue field with characteristic not equal to two. 
\subsection*{Acknowledgements} 
We are very grateful to Michael Joswig, Kristan Ranestad and Bernd Sturmfels for many interesting discussions on the project and valuable comments on earlier versions of this paper.

\section{Tropical surfaces and tropical lines}\label{section:trop}

In this section  we go through the basic definitions and we fix  our notation concerning tropical surfaces and tropical lines in $\rr^3$. We refer to \cite{MaclaganSturmfels15} for further reading on the topic. We work over the {\em tropical semiring} $(\rr \cup \{\infty\},\oplus, \odot)$, where the tropical operations are defined as  $x \oplus y= \min\{x,y\}$ and $x\odot y = x+y$. 
Consider a tropical polynomial  
\[f(x,  y, z) \quad = \!\! \bigoplus_{\mathbf{v}=(v_1, v_2,v_3) \in \zz^3} \!\!\!\!\! c_{\mathbf{v}} \odot
 x_1^{\odot v_1} \odot x_2^{\odot v_2} \odot x_3^{\odot v_3} \,\,\,= \,\,\,\bigoplus_{\mathbf{v} \in \zz^3} c_{\mathbf{v}} x^{\mathbf{v}}, \]
 with $c_v \not = \infty$  for finitely many coefficients.  The \emph{tropical surface $T(f)$} is defined as the set of points in $\rr^3$ at which the minimum  among the quantities $c_{\mathbf{v}} + \mathbf{v} \cdot x$ is attained at least twice. It is well known that $T(f)$ is a pure connected polyhedral complex of dimension $2$, some of whose cells are unbounded in $\rr^3$.

The \emph{Newton polytope of $f$} is the lattice polytope 
\[ \text{Newt}(f) =  \text{conv} \big\{\mathbf{v} \, : \, c_{\mathbf{v}} \not = \infty \big\}. \]
Let $\mathcal{A} = \text{Newt}(f) \cap \zz^3$ be its set of  lattice points. The coefficients of $f$ induce a \emph{regular subdivision} $\ks$ of $\mathcal{A}$ by taking the convex hull in $\rr^{4}$ of the points $(\mathbf{v}, c_{\mathbf{v}})$ and projecting the lower faces to $\text{Newt}(f)$.
The coefficient vectors inducing the same subdivision $\ks$ form a relatively open polyhedral cone $\Sigma(\ks)$ in $\rr^{|\mathcal{A}|}$, called \emph{the secondary cone}, see \cite[Chapter 5]{DRS10}.  The tropical hypersurface $T(f)$ is dual to the subdivision $\ks$,  and they determine each other \cite[Proposition 3.1.6]{MaclaganSturmfels15}. More precisely, there is an inclusion-reversing bijection between the polyhedral cells in $T(f)$ and the cells in the subdivision $\ks$: $k$-dimensional cells of $T(f)$ are mapped to $3-k$ dimensional cells in $\ks$.  We say that $T(f)$ is \emph{smooth} if the subdivision is a unimodular triangulation, i.e.,   all maximal cells are tetrahedra of lattice  volume one. For basics definitions on triangulations we refer to \cite{DRS10}. 

Given $d \in \nn$, the polytope $d \, \Delta_3$ is  $\text{conv}\big( (0,0,0), (d,0,0),(0,d,0),(0,0,d) \big)$. By a subdivision of $d\, \Delta_3$ we will always mean a subdivision of its lattice points. The vectors $-(e_1+e_2+e_3)$, $e_1$, $e_2$ and $e_3$ are denoted $\omega_0$, $\omega_1$, $\omega_2$ and $\omega_3$.  For $i\in \{0,\dotsc,3\}$, we use the notation $F_i$ for  the facet of $d\, \Delta_3$ with inner normal vector $\omega_i$.  Moreover, for distinct $i,j\in\{0,\dotsc,3\}$, we set $F_{ij}:=F_i\cap F_j$. 

 A {\em tropical surface of degree $d$} is a subset of $\rr^3$ of the form $T(f)$, where $f$ is a tropical polynomial whose Newton polytope is $d\,\Delta_3$, see \cite[Chapter 4.5]{MaclaganSturmfels15}. If $X = T(f)$ is a tropical surface of degree $d$, then we indicate with $\ks_X$ the regular subdivision induced by the coefficients of $f$.  Let $\ks$ be regular subdivision of $d\,\Delta_3$ for some $d\in\nn$, and $\Sigma(\ks)$ its secondary cone. If $c$ is a point in $\Sigma(\ks)$, then $X_c$ is the associated tropical surface. More precisely, $c = (c_{{\bf v}})_{{\bf v} \in d\, \Delta_3}$ and  $X_c:=T(f)$, where 
\[f\quad = \,\,\,\bigoplus_{\mathbf{v} \in d\,\Delta_3} c_{\mathbf{v}} x^{\mathbf{v}}. \]

Let $L$ be an unrooted tree with five edges, and six vertices, two of which are 3-valent and the rest 1-valent. We define a {\em tropical line in $\rr^3$} to be any realization of $L$ in $\rr^3$ such that the realization is a balanced polyhedral complex with four unbounded rays (the 1-valent vertices of $L$ are pushed to infinity) having direction vectors  $\omega_0$, $\omega_1$, $\omega_2$ and $\omega_3$. Balanced means that  at each vertex the primitive integer vectors in the directions of all outgoing edges adjacent to it  sum to zero. If the bounded edge has length zero, the tropical line is called {\em degenerate}. 

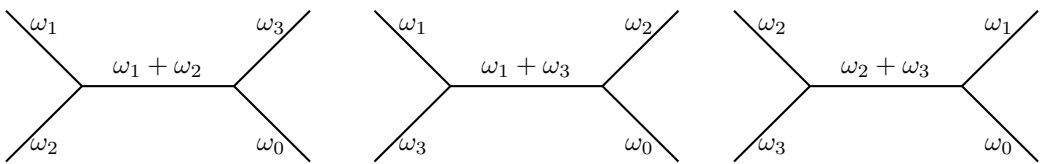
\begin{figure}[h]
\begin{center}
\begin{tikzpicture}
\draw [thick] (-1,1)--(0,0);
\draw [thick] (-1,-1)--(0,0);
\draw [thick] (0,0) --(2,0);
\draw [thick] (2,0) --(3,1);
\draw [thick] (2,0) --(3,-1);
\node [above] at (1,0) {$\omega_1+\omega_2$};
\node [right] at (-0.8,0.8) {$\omega_1$};
\node [right] at (-0.8,-0.8) {$\omega_2$};
\node [left] at (2.8,0.8) {$\omega_3$};
\node [left] at (2.8,-0.8) {$\omega_0$};  
\end{tikzpicture}
\ \ \ \ \
\newcommand{\drawskeleton}{
\draw [thick] (-1,1)--(0,0);
\draw [thick] (-1,-1)--(0,0);
\draw [thick] (0,0) --(2,0);
\draw [thick] (2,0) --(3,1);
\draw [thick] (2,0) --(3,-1);
}
\begin{tikzpicture}
\drawskeleton
\node [above] at (1,0) {$\omega_1+\omega_3$};
\node [right] at (-0.8,0.8) {$\omega_1$};
\node [right] at (-0.8,-0.8) {$\omega_3$};
\node [left] at (2.8,0.8) {$\omega_2$};
\node [left] at (2.8,-0.8) {$\omega_0$};  
\end{tikzpicture}
\ \ \ \ \ 
\begin{tikzpicture}
\draw [thick] (-1,1)--(0,0);
\draw [thick] (-1,-1)--(0,0);
\draw [thick] (0,0) --(2,0);
\draw [thick] (2,0) --(3,1);
\draw [thick] (2,0) --(3,-1);
\node [above] at (1,0) {$\omega_2+\omega_3$};
\node [right] at (-0.8,0.8) {$\omega_2$};
\node [right] at (-0.8,-0.8) {$\omega_3$};
\node [left] at (2.8,0.8) {$\omega_1$};
\node [left] at (2.8,-0.8) {$\omega_0$};  
\end{tikzpicture}
\caption{The types of non-degenerate tropical lines in $\rr^3$.}\label{spacelines} 
\end{center}
\end{figure}

For non-degenerate tropical lines, there are three types of tropical lines in $\rr^3$, as shown in Figure \ref{spacelines}. The  types of the lines in Figure \ref{spacelines}, from left to right, are denoted by $((12)(30))$, $((13)(20))$ and $((23)(01))$. Each innermost pair of digits indicate the directions of two adjacent rays.

In classical geometry, any two distinct points lie on a unique line. When we turn to tropical lines, this is true only for generic points. In fact, for special choices of points $P$ and $Q$ there are infinitely many tropical lines passing through $P$ and $Q$. We call an  infinite collection of tropical lines in $\rr^3$  a {\em two-point family} if there exist two points lying on all lines in the collection.

\begin{lemma}\label{lem:uniqueline}
Let $P,Q\in\rr^3$. There exist infinitely many tropical lines containing $P$ and $Q$ if and only if one of the coordinates of the vector $Q-P$ is zero, or two of them coincide. In all other cases, $P$ and $Q$ lie on a unique line.
\end{lemma}

The group $S_4$ of permutations of four elements  acts naturally on many of the spaces involved with tropical surfaces. Firstly, observe that $S_4$ is the symmetry group of the simplex $d\,\Delta_3$. This induces an action of the set of lattice points $d\,\Delta_3 \cap \zz^3$ (in fact on all of $\zz^3$).  Obviously, this action of $S_4$ also induce an action of $S_4$ on the set of subdivisions of $d\,\Delta_3$.

Secondly, $S_4$ acts on the set of tropical surfaces of degree $d$. Let $X=T(f)$, where $f(x)= \bigoplus_{\mathbf{v} \in d\, \Delta_3} c_{\mathbf{v}} x^{\mathbf{v}}$. For a given $\sigma\in S_4$, we define $\sigma(X)$ to be the surface $T(\sigma(f))$, where  $\sigma(f)(x)= \bigoplus_{\mathbf{v} \in d\,\Delta_3} c_{\sigma(\mathbf{v})} x^{\sigma(\mathbf{v})}$. 
Clearly, $\sigma(X)$ is still of degree $d$, and the resulting action is compatible with the action of $S_4$ on the subdivisions of $d\,\Delta_3$. In other symbols, $\ks_{\sigma(X)}=\sigma(\ks_X)$.

\section{Motifs of tropical lines}  \label{sec:lsc}
Before focusing on tropical cubic surfaces, we begin by studying tropical lines contained in tropical surfaces of arbitrary degree. 
It is important to note that containment is meant purely set-theoretically. The symbols $X$ and $L$ will always refer to  a tropical surface of degree $d$ in $\rr^3$, and a tropical line, respectively. Hence in particular, the statement $L\sub X$ means that $L$ is contained in $X$ as subsets of $\rr^3$. 

Given a cell $C$ of the tropical surface $X$, we use the notation $C^\vee$ to indicate the corresponding dual cell in the subdivision $\ks_X$, and vice versa. More generally, for any tropical surface $X$, we have the natural map $\mc_X\colon |X| \pil X$, taking a point $p$ in the support of $X$ to the minimal cell of $X$ containing $p$. Combining $\mc_X$ with dualization, we get the map 
\begin{equation}
\begin{split}
\mc_X^\vee\colon |X|& \longrightarrow \ks_X\\
p& \longmapsto \mc_X(p)^\vee.
\end{split}
\end{equation}
If $Y\sub X$ is any subset, we set $\mc_X^\vee(Y):=\bigcup_{y\in Y}\mc_X^\vee(y).$
Note that if $Y$ is connected, then $\mc_X^\vee(Y)$ is a connected subcomplex of $\ks_X$. Whenever the surface $X$ is clear from the context, we will omit the symbol $X$ from the maps $\mc_X$ and $\mc_X^\vee$. 
\vspace{\baselineskip}

The notion of {\em trespassing line segments}  on $X$ was introduced in \cite{Vigeland10} by the second author. If $\ell\sub X$ is any ray or line segment, we say that $\ell$ is trespassing on $X$ if there exist distinct cells $C_1,C_2 \sub X$ such that 
\begin{equation}\label{eq:tres}
\dim\big(\inter(C_1) \cap \ell\big)=\dim\big(\inter(C_2) \cap \ell\big)=1.
\end{equation} 
Alternatively, $\ell$ is trespassing on $X$ if it is not contained in the closure of any single cell of $X$. For smooth $X$, trespassing can happen in one way only, as shown in \cite[Lemma 6.2]{Vigeland10}:

\begin{lemma}\label{lem:tres}
Let $\ell \sub X$, where $\ell$ is a line segment and $X$ is smooth. If $C_1,C_2\sub X$ are any cells satisfying \eqref{eq:tres} and such that $\ell\sub C_1\cup C_2$, then $C_1$ and $C_2$ are 2-dimensional cells of $X$ whose intersection is a vertex $V$ of $X$. 
\end{lemma}  
An immediate consequence of this is that $C_1^\vee$ and $C_2^\vee$ are opposite edges of the tetrahedron $V^\vee$ in the unimodular triangulation $\ks_X$. The following converse to Lemma \ref{lem:tres} is straightforward:

\begin{lemma}\label{lem:tresconv}
Let $e$ and $e'$ be opposite edges of a tetrahedron $\tau$ in the unimodular triangulation  $\ks_X$. Then there is a trespassing line segment on $X$ passing through the vertex $\tau^\vee\in X$, and which is orthogonal to both $e$ and $e'$.
\end{lemma}

We recall one more result  \cite[Lemma 6.4]{Vigeland10}:
\begin{lemma}\label{lem:1vert}
Suppose $L\sub X$ is non-degenerate, and that the vertex $v$ of $L$ lies in the interior of a 1-dimensional cell $E$ of $X$. Then $L\cap E=\{v\}$, and the three edges of $L$ adjacent to $v$ start off in different 2-dimensional cells of $X$ adjacent to $E$. 
%
\end{lemma}
Equivalently, suppose that $\omega_i$ and $\omega_j$ are the direction vectors of the unbounded edges of $L$ emanating from $v$. Then the edges of the triangle $E^\vee$ are orthogonal to $\omega_i$, $\omega_j$ and $\omega_i+\omega_j$ respectively.

\vspace{\baselineskip}

We next describe a way of displaying the essential information of how $L$ lies on $X$. 
For any tropical line $L\sub \rr^3$, a {\em decoration}  consists of a finite number of dots (possibly none) on each edge, and at each vertex either a dot, a vertical line segment, or nothing. A {\em marked tropical line $L$} is a tropical line with decorations applied as follows: 
\begin{enumerate}[(a)]
\item If an edge of $L$ is trespassing, we marked the edge with a dot at each point  where the edge is trespassing. 
\item For each vertex $v$ of $L$, the corresponding vertex of the graph has a dot if $\dim(\mc_X(v))=0$, and a vertical line segment if $\dim(\mc_X(v))=1$.    
\end{enumerate}
The {\sl combinatorial type} of a marked tropical line $L$ is a decorated graph $G$ encoding the directions of the unbounded edges and the decoration of $L$.  We consider the graphs without metrics, so moving an edge-dot along the interior of its  edge does not change the decoration.

\begin{definition}
Let $X$ be a smooth tropical surface, and let $L\sub X$ be any marked tropical line. The {\em motif $\mathcal{M}$\, of $L$ on $X$} is given by a pair  $\mathcal{M} = (G,\mathcal{R})$, where $G$ is the combinatorial type of the marked line $L$ and $\mathcal{R}$ is the subcomplex of $\ks_X$ such that  $\mc_X^\vee(L)=\mathcal{R}$. We refer to $G$ as the primal motif and to $\mathcal{R}$ as the dual motif. 

Conversely, suppose $\mathcal{M} = (G, \mathcal{R})$ is a motif with $\mathcal{R}$ a subcomplex of a unimodular regular triangulation~$\mathcal{S}$. If $X'$ is any tropical surface with $\ks_{X'}=\ks$, we say that $\mathcal{M}$ is {\em realized on $X'$} if there is a marked tropical line $L\sub X'$ with combinatorial type $G$ such that $\mc_{X'}^\vee(L)=\mathcal{R}$.  
\end{definition} 
We will often consider primal motifs  modulo automorphisms of decorated graphs. 

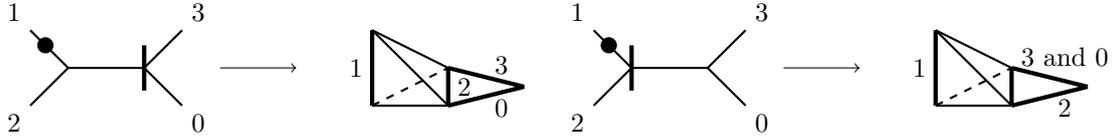
\begin{figure}[ht]
\begin{center}
\begin{tikzpicture}
\draw [thick] (-0.5,0.5)--(0,0);
\draw [thick] (-0.5,-0.5)--(0,0);
\draw [thick] (0,0) --(1,0);
\draw [thick] (1,0) --(1.5,0.5);
\draw [thick] (1,0) --(1.5,-0.5);
\draw [ultra thick] (1,0.3)--(1,-0.3);
\draw [fill] (-0.3,0.3) circle [radius=.1];
 \node [above left] at (-0.5,0.5) {$1$};
 \node [below left] at (-0.5,-0.5) {$2$};
 \node [above right] at (1.5,0.5) {$3$};
 \node [below right] at (1.5,-0.5) {$0$};
 \draw [->] (2,0) -- (3,0);
 \draw [ultra thick] (4,0.5)--(4,-0.5);
 \draw [thick] (4,-0.5)--(5,-0.5);
  \draw [thick, dashed] (4,-0.5)--(5,0);
  \draw [ultra thick] (5,0)--(5,-0.5);
   \draw [thick] (4,0.5)--(5,-0.5);
  \draw [thick] (4,0.5)--(5,0);
  \draw [ultra thick] (5,0)--(6,-0.25);
   \draw [ultra thick] (5,-0.5)--(6,-0.25);
    \node [left] at (4,0) {$1$};
 \node [right] at (5,-0.25) {$2$};
 \node [above] at (5.7,-0.2) {$3$};
 \node [below] at (5.7,-0.3) {$0$};
\end{tikzpicture}
\,\,\,\,
\begin{tikzpicture}
\draw [thick] (-0.5,0.5)--(0,0);
\draw [thick] (-0.5,-0.5)--(0,0);
\draw [thick] (0,0) --(1,0);
\draw [thick] (1,0) --(1.5,0.5);
\draw [thick] (1,0) --(1.5,-0.5);
\draw [ultra thick] (0,0.3)--(0,-0.3);
\draw [fill] (-0.3,0.3) circle [radius=.1];
 \node [above left] at (-0.5,0.5) {$1$};
 \node [below left] at (-0.5,-0.5) {$2$};
 \node [above right] at (1.5,0.5) {$3$};
 \node [below right] at (1.5,-0.5) {$0$};
 \draw [->] (2,0) -- (3,0);
 \draw [ultra thick] (4,0.5)--(4,-0.5);
 \draw [thick] (4,-0.5)--(5,-0.5);
  \draw [thick, dashed] (4,-0.5)--(5,0);
  \draw [ultra thick] (5,0)--(5,-0.5);
   \draw [thick] (4,0.5)--(5,-0.5);
  \draw [thick] (4,0.5)--(5,0);
  \draw [ultra thick] (5,0)--(6,-0.25);
   \draw [ultra thick] (5,-0.5)--(6,-0.25);
    \node [left] at (4,0) {$1$};
 \node [above] at (5.7,-0.1) {$3$ and $0$};
 \node [below] at (5.7,-0.3) {$2$};
\end{tikzpicture}
\caption{Two motifs giving the same dual motif structure, but with different sets of required exits (indicated by bold lines). The exits are determined using Lemmas \ref{lem:tres} and \ref{lem:1vert}.}\label{fig:lscex}
\end{center}\end{figure}

\medskip 

Because tropical lines in $\rr^3$ are unbounded, any dual motif in $\ks_X$  contains cells dual to unbounded cells of $X$. Such cells of $\ks_X$ are precisely those lying in the boundary of $d\,\Delta_3$. Let $\tau$ be a lattice polytope (of dimension 1,2 or 3) contained in $d\,  \Delta_3$. We say that $\tau$ has an {\em exit in the direction of $\omega_i$} if at least one edge of $\tau$ lies in $F_i$. If $\tau$ has exits in the directions of $m$ of the $\omega_i$'s, we say that $\tau$ has {\em $m$ exits}.

The relevance of this definition should be clear from the following observation. 
\begin{remark}
 Let $C$ be any cell of $X$, and let $p\in C$ be an arbitrary point. Then $C$ contains the ray with starting point $p$ and direction $\omega_i$ if and only if $C^\vee$ has an exit in direction $\omega_i$. 
 \end{remark}

When $X$ is smooth, the cell structure of a  dual motif $\mc_X^\vee(L)$ is in many cases uniquely determined by the primal motif $G$ of $L$ on $X$. Moreover, using Lemma \ref{lem:tres} and Lemma \ref{lem:1vert}, we can often describe explicitly the exits required of the edges of $\mc_X^\vee(L)$. For example, the two  motifs in Figure  \ref{fig:lscex} imply the same cell structure of $\mc_X^\vee(L)$, but with different exit properties.

\begin{remark}
A dual motif $\mathcal{R}\subset \ks$ will often have more exits than those required by the position. Hence it is usually more difficult to reverse the process described in the last paragraph, i.e., to determine the primal motif of $L$ on $X$, given  $\mc_X^\vee(L)\sub\ks_X$.
\end{remark}

\begin{remark} \label{rem:middleedge}
It is important to mention that there is  one case where the cell structure of $\mc_X^\vee(L)$ is not determined by the primal motif of $L$ on $X$. Namely, when both vertices of $L$ are vertices of $X$, and the middle edge of $L$ is not trespassing. In this case, the middle edge of $L$ may or may not be an edge of $X$, giving different structures of $\mc_X^\vee(L)$, see Figure \ref{fig:tvecp}.

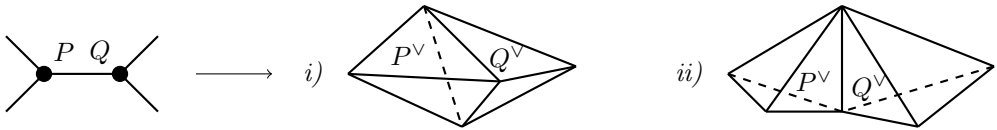
\begin{figure}[ht]
\begin{center}
\begin{tikzpicture}
\draw [thick] (-0.5,0.5)--(0,0);
\draw [thick] (-0.5,-0.5)--(0,0);
\draw [thick] (0,0) --(1,0);
\draw [thick] (1,0) --(1.5,0.5);
\draw [thick] (1,0) --(1.5,-0.5);
\draw [fill] (0,0) circle [radius=.1];
\draw [fill] (1,0) circle [radius=.1];
 \node [above right] at (0,0.05) {$P$};
 \node [above left] at (1,0) {$Q$};
 \draw [->] (2,0) -- (3,0);
  \node [left] at (3.8,0) {\emph{i)}};
 \node [above] at (4.8,0) {$P^\vee$};
 \node [above] at (6.1,-0.1) {$Q^\vee$};
 \draw [thick] (4,0)--(6,-0.1);
 \draw [thick] (4,0)--(5.5,-0.7);
 \draw [thick](5.5,-0.7)--(6,-0.1);
 \draw [thick]((4,0)--(5,0.9);
 \draw [thick](5,0.9)--(6,-0.1);
 \draw [thick] (6,-0.1)--(7,0.1);
 \draw [thick,dashed](5,0.9)--(5.5,-0.7);
 \draw[thick](5,0.9)--(7,0.1);
 \draw[thick](5.5,-0.7)--(7,0.1);
  \node [left] at (8.8,0) {\emph{ii)}};
 \draw[thick](9,0)--(9.5,-0.5);
 \draw[thick](9.5,-0.5)--(10.5,-0.5);
 \draw[thick](10.5,-0.5)--(10.5,0.9);
 \node [left] at (10.5,-0.1) {$P^\vee$};
 \node [below right ] at (10.5,0.1) {$Q^\vee$};
 \draw[thick](9,0)--(10.5,0.9);
 \draw[thick](9.5,-0.5)--(10.5,0.9);
 \draw[thick](10.5,-0.5)--(11.5,-0.7);
 \draw[thick](12.5,0.1)--(11.5,-0.7);
\draw[thick](12.5,0.1)--(10.5,0.9);
 \draw[thick](11.5,-0.7)--(10.5,0.9);
 \draw[thick,dashed](12.5,0.1)--(10.5,-0.5);
 \draw[thick,dashed](9,0)--(10.5,-0.5);
\end{tikzpicture}
\caption{A primal motif giving two possible dual motifs, depending on whether $PQ$ is a 1-dimensional cell of $X$ (case $i$) or not (case $ii$).}\label{fig:tvecp}
\end{center}\end{figure}
The two tetrahedra $P^\vee$ and $Q^\vee$ have a common facet if $PQ$ is an edge of $X$ (case \emph{i}), but only an edge in common otherwise (i.e., if $PQ$ goes across a 2-dimensional cell of $X$).
Note that if the middle edge were trespassing, there would be no ambiguity: By Lemma \ref{lem:tres}, no point of $PQ$ could then be in the interior of a 1-dimensional cell of $X$.
\end{remark}
\vspace{\baselineskip}

We finish this section by introducing the definitions of deformations and specializations of tropical lines on a surface. 
Let $\ks$ be a given subdivision of $d \, \Delta_3$, and let $\Sigma(\ks)$ be the corresponding secondary cone. We define the incidence complex $\xx_\ks\sub \Sigma(\ks)\times \rr^3$ by
\begin{equation*}
  \xx_\ks:=\{(c,x)\:|\:x\in X_c\}\sub \Sigma(\ks)\times \rr^3.
\end{equation*}
Using the Euclidean metric on both $\Sigma(\ks)$ and $\rr^3$, we give $\xx_\ks$ the topology induced by the product topology on $\Sigma(\ks)\times \rr^3$. This makes the projections on $\Sigma(\ks)$ and $\rr^3$, denoted by $p_1$ and $p_2$ respectively, continuous. Note that for any $c\in \Sigma(\ks)$, the set $p_2(p_1^{-1}(c))$ is the tropical surface $X_c$.

\begin{definition}
A {\em family of tropical lines associated to $\ks$} is a subset $\LL\sub \xx_\ks$ satisfying the following conditions:
\begin{itemize}
\item For any $c \in \Sigma(\ks)$, $p_2(p_1^{-1}(c)\cap \LL)$ is a tropical line $L_c \sub X_c$.
\item The projections from $\LL$ to $\Sigma(\ks)$ and $\rr^3$ are continuous.
\end{itemize}
\end{definition}

\begin{definition}
A {\em deformation} of $L\sub X_c$ is a family $\LL$ of tropical lines associated to $\ks$, such that 
\begin{itemize}
\item $p_1(\LL)$ contains $\alpha$, and is homeomorphic to an interval,
\item for any two points $c_1\neq c_2$ in $p_1(\LL)$, we have $X_{c_1}\neq X_{c_2}$.
\end{itemize}
\end{definition}
Note that a deformation of $L\sub X$ can be thought of as a map $t\mapsto (L_t,X_t)$, where $t$ runs through some interval $I\sub \rr$ containing $0$, and where $(L,X)=(L_0,X_0)$. In particular, $0$ can be an endpoint of $I$, as in $I=[0,1)$.

\begin{definition}
Let $L$ be a tropical line with motif $\cc$ on $X$. We say that $L$ {\em deforms into motif $\cc'$}, if there exist a deformation $t\mapsto (L_t,X_t)$ of $L\sub X$ such that for all $t\in I\smallsetminus \{0\}$, the motif of $L_t$ on $X_t$ is $\cc'$.
\end{definition}

The following lemma gives a simple property of deformations, namely that one cannot deform a tropical line away from
a vertex through which it is trespassing.
\begin{lemma}\label{lem:deftres}
Let $X$ be smooth. Suppose $L\sub X$ has a trespassing line segment $\ell$, passing through a vertex $V$ which is dual to the tetrahedron  $\Delta$ in $\ks_X$. For any deformation $t\mapsto (L_t,X_t)$, $t\in I$ of $L\sub X$, let $\ell_t$ be the edge of $L_t$ parallel to $\ell$. Then for $t$ small enough, $\ell_t$ is trespassing through $\Delta^\vee\in X_t$.
\end{lemma}
\begin{proof}
Since $\ell$ is trespassing through $V$, Lemma \ref{lem:tres} gives that $\dim(\ell \, \cap \,\inter(\La^\vee)) = \dim(\ell \, \cap\, \inter(\La'^\vee))=1$, for some pair of opposite edges $\La,\La'$ of $\Delta$. By continuity of the deformation, this implies that $\dim(\ell_t\cap\inter(\La^\vee))=\dim(\ell_t\cap\inter(\La'^\vee))=1$ for small enough $t$. Hence $\ell_t$ is trespassing through $\Delta^\vee$ in $X_t$. 
\end{proof}

\begin{remark}
Note that the proof of Lemma \ref{lem:deftres} rests on Lemma \ref{lem:tres}, which requires $X$ to be smooth. In fact, it is possible to produce examples of non-smooth $X$ where one {\em can} deform away from trespassed vertices.
\end{remark}

Related to deformations is the concept of specialization:
\begin{definition}
Let $t\mapsto (L_t,X_t)$ be a deformation of $L_0\sub X_0$, where $t\in [0,1]$. We say that $L_0\sub X_0$ {\em specializes} to $L_1\sub X_1$ if the motif of $L_t\sub X_t$ is constant for all $t\in[0,1)$ but differs for $t=1$.
\end{definition}
We will use these notions and results in Section \ref{sec:cubic} where we will analyze the surfaces dual to a honeycomb triangulation. 

\section{Preliminary results} 

Let $d\in\nn$ be fixed, and let $\ks$ be a  regular unimodular triangulation of $d\,\Delta_3$.  Recall that each $c \in \Sigma(\ks)$ corresponds to a smooth tropical surface $X_c$ with subdivision $\ks$. 
We say that a property $\mathcal{P}$ holds for {\em general smooth tropical surfaces with subdivision $\ks$} if $\mathcal{P}$ holds for $X_c$ for every $c$ in some open dense subset of the secondary cone $\Sigma(\ks)$. 
More generally, $\mathcal{P}$ holds for {\em general smooth tropical surfaces of degree $d$} if $\mathcal{P}$ holds for general tropical surfaces with triangulation $\ks$, for all regular unimodular  triangulations $\ks$ of $d\,\Delta_3$.
Finally, $\mathcal{P}$ holds for {\em general smooth tropical surfaces} if $\mathcal{P}$ holds for general smooth tropical surfaces of degree $d$, for all $d\in\nn$.

The next lemma gives an important example of a property held by general smooth tropical surfaces (in all degrees).

\begin{lemma}\label{lem:doublytres}
A general smooth tropical surface $X$ contains no doubly trespassing line segments.
\end{lemma}
\begin{proof}
Let $X=T(f)$ be a smooth tropical surface, given by a tropical polynomial $f=\bigoplus c_{{\bf v}} x^{{\bf v}}$, and suppose $\ell\sub X$ is a line segment containing two vertices of $X$, say $P$ and $Q$, in its relative interior. We will show that this condition implies that there is a linear relation between the coefficients $c_{{\bf v}}$.

The situation is shown in Figure \ref{fig:doublytres}. From Lemma \ref{lem:tres}, it follows that $\ks_X$ contains three $1$-dimensional cells $AB,CD,EF$ such that $\ell\sub (AB)^\vee\cup(CD)^\vee\cup (EF)^\vee$, and such that $P^\vee=ABCD$ and $Q^\vee=CDEF$. Obviously, a necessary condition for this to happen is that $\ell$ is parallel to both vectors products $\overrightarrow{AB}\times \overrightarrow{CD}$ and $\overrightarrow{CD}\times \overrightarrow{EF}$, see Lemma \ref{lem:tresconv}.
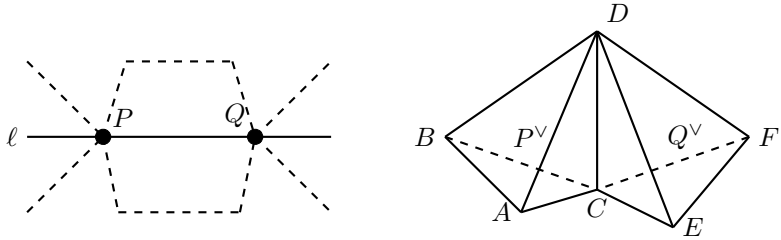
\begin{figure}[h]
\begin{tikzpicture}
\draw [thick] (-1,0)--(3,0);
 \draw [thick,dashed](-1,1)--(0,0);
 \draw [thick,dashed](-1,-1)--(0,0);
  \draw [thick,dashed](0,0)--(0.3,1);
 \draw [thick,dashed](0.3,1)--(1.7,1);
 \draw [thick,dashed](1.7,1)--(2,0);   
   \draw [thick,dashed](0,0)--(0.2,-1);
 \draw [thick,dashed](0.2,-1)--(1.8,-1);
 \draw [thick,dashed](1.8,-1)--(2,0);   
 \draw [thick,dashed](2,0)--(3,1);
 \draw [thick,dashed](2,0)--(3,-1);
 \draw [fill] (0,0) circle [radius=.1];
\draw [fill] (2,0) circle [radius=.1];
\node [above right] at (0,0) {$P$};
\node [above left] at (2,0) {$Q$};
\node [left] at (-1,0) {$\ell$};

 \draw[thick](4.5,0)--(5.5,-1);
 \draw[thick](5.5,-1)--(6.5,-0.7);
 \draw[thick](6.5,-0.7)--(6.5,1.4);
 \node [left] at (6,0) {$P^\vee$};
 \node [right ] at (7.3,0) {$Q^\vee$};
 \draw[thick](4.5,0)--(6.5,1.4);
 \draw[thick](5.5,-1)--(6.5,1.4);
 \draw[thick](6.5,-0.7)--(7.5,-1.2);
 \draw[thick](8.5,0)--(7.5,-1.2);
\draw[thick](8.5,0)--(6.5,1.4);
 \draw[thick](7.5,-1.2)--(6.5,1.4);
 \draw[thick,dashed](8.5,0)--(6.5,-0.7);
 \draw[thick,dashed](4.5,0)--(6.5,-0.7);
 \node [left ] at (5.5,-1) {$A$};
 \node [left ] at (4.5,0) {$B$};
 \node [below ] at (6.5,-0.7) {$C$};
 \node [above right ] at (6.5,1.4) {$D$};
  \node [right ] at  (7.5,-1.2) {$E$};
  \node [right ] at (8.5,0) {$F$};
       
\end{tikzpicture}

\caption{A doubly trespassing line, and the dual configuration.}
\label{fig:doublytres}
\end{figure}

Since $\ell$ is contained in each of the planes spanned by $(AB)^\vee$, $(CD)^\vee$ and $(EF)^\vee$, any point $p\in\ell$ satisfies the defining equations of these planes: 
\begin{eqnarray}\label{eqs:la}
c_{A}+A\cdot x&=c_{B}+B \cdot x\\
c_{C}+C \cdot x&=c_{D}+D \cdot x\\
c_{E}+ E \cdot x&=c_{F}+B\cdot x.
\end{eqnarray}
Substituting $p$ for $x$ and rearranging, this amounts to the matrix equation
\begin{equation}\label{eq:MLP}
M\begin{pmatrix}
p\\
1
\end{pmatrix}=0,\quad\text{where } 
M=
\begin{pmatrix}
\vec{AB}\;&c_B-c_A\\
\vec{CD}\;&c_D-c_C\\
\vec{EF}\;&c_F-c_E    
\end{pmatrix}=
\begin{pmatrix}
B_1-A_1&B_2-A_2&B_3-A_3&c_B-c_A\\
D_1-C_1&D_2-C_2&D_3-C_3&c_D-c_C\\
F_1-E_1&F_2-E_2&F_3-E_3&c_F-c_E    
\end{pmatrix}.
\end{equation} 
Since \eqref{eq:MLP} holds for all $p$ in $\ell$, the nullity of $M$ must be at least 2, i.e., $\rank M\leq 2$. In fact, $\rank M=2$, since the vectors $\overrightarrow{AB}$, $\overrightarrow{CD}$, $\overrightarrow{EF}$ are not parallel. Therefore there is a linear relation between $c_A,\dotsc,c_F$.

To prove the lemma, let $\ks$ be an  unimodular triangulation $\ks$ of $d\, \Delta_3$. In $\ks$, we look for all pairs of tetrahedra with a common edge, and such that (with the notation of Figure \ref{fig:doublytres}) $\overrightarrow{AB}\times \overrightarrow{CD}$ is parallel to $\overrightarrow{CD}\times \overrightarrow{EF}$. There are at most finitely many such pairs. As seen above, each such pair gives rise to a hyperplane section of $\Sigma(\ks)$, and any surface containing a doubly trespassing line, corresponds to a point on one of these hyperplanes. This proves the lemma.
\end{proof}

The above lemma greatly limits the number of ways in which a tropical line can lie on a general smooth tropical surface. In particular, the lemma says that for general $X$, each of the five edges of $L\sub X$ contains at most one vertex of $X$ in its relative interior. An immediate implication of this is the following interesting result: There exists a {\em finite} list encompassing motifs of  lines on general smooth surfaces of any degree.

There are some motifs of $L$ that do not occur on general $X$, but which are not excluded by Lemma \ref{lem:doublytres}. Many of these can be identified using the lemma to follow.  By a {\em 3-star on $X$} we will mean the union of 3 line segments, pairwise non-parallel, with a common endpoint, all of which is pointwise contained in $X$.

\begin{lemma}\label{lem:special}
A general smooth tropical surface $X$ in $\mathbb{R}^3$ does not contain the 3-stars in Figure \ref{fig:closedcases}. The vertices $V_1$ and $V_2$ in case  iii) are not adjacent on $X$. 
\end{lemma}

\begin{figure}[h]
\centering

\begin{tikzpicture}
 \node [left] at (-1,0) {\emph{i)}};
\draw [thick] (-1,1)--(0,0);
\draw [thick] (-1,-1)--(0,0);
\draw [thick] (0,0) --(2,0);
\draw [fill] (-0.5,0.5) circle [radius=.1];
\draw [fill] (-0.5,-0.5) circle [radius=.1];
\draw [fill] (0.5,0) circle [radius=.1];
 \node [above right] at (0,0) {$v$};
   \end{tikzpicture}
   \,\,\,
 \begin{tikzpicture}
  \node [left] at (-1,0) {\emph{ii)}};
\draw [thick] (-1,1)--(0,0);
\draw [thick] (-1,-1)--(0,0);
\draw [thick] (0,0) --(2,0);
\draw [fill] (0,0) circle [radius=.1];
\draw [fill] (1,0) circle [radius=.1];
 \node [above right] at (0,0) {$v$};
   \end{tikzpicture}
   \,\,\,
\begin{tikzpicture}
  \node [left] at (-1.2,0) {\emph{iii)}};
\draw [thick] (-1,1)--(0,0);
\draw [thick] (-1,-1)--(0,0);
\draw [thick] (0,0) --(2,0);
 \node [ below left] at (-0.5,0.5) {$V_1$};
  \node [ above left] at (-0.5,-0.5) {$V_2$};
\draw [fill] (-0.5,0.5) circle [radius=.1];
\draw [fill] (-0.5,-0.5) circle [radius=.1];
\draw [ultra thick] (0,0.5)--(0,-0.5);
 \node [above right] at (0,0) {$v$};
 \end{tikzpicture}
\caption{The 3-stars considered in Lemma \ref{lem:special}.} \label{fig:closedcases}
\end{figure}
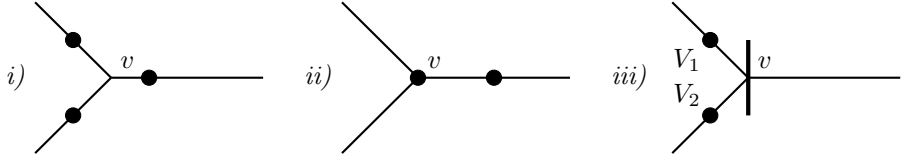
\begin{proof}
Let $X=T(f)$ be smooth of some arbitrary degree $d$, where $f= \oplus_{{\bf v}} c_{{\bf v}} x^{{\bf v}}$. For any 3-star $Y\sub X$ in Figure  \ref{fig:closedcases}, we  consider the subcomplex $\mc^\vee(Y)$ in $\ks_X$ depicted in  Figure \ref{fig:closedsubs}. 

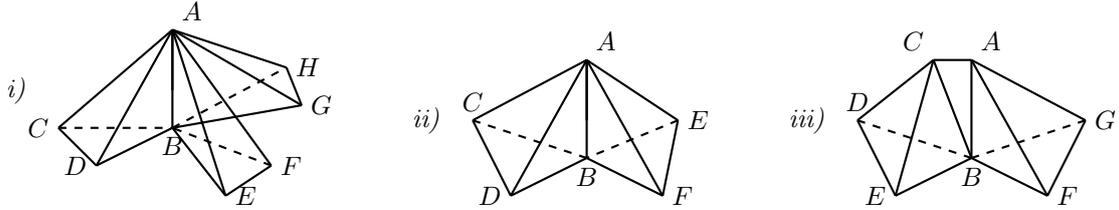
\begin{figure}[h]
\begin{tikzpicture}
 \node [left] at  (0.7,0) {\emph{i)}};
 \draw[thick](1,-0.5)--(1.5,-1);
 \draw[thick](1.5,-1)--(2.5,-0.5);
 \draw[thick](2.5,-0.5)--(2.5,0.8);
 \draw[thick](2.5,0)--(2.5,0.8);
 \draw[thick](1.5,-1)--(2.5,0.8);
 \draw[thick](2.5,-0.5)--(3.2,-1.4);
 \draw[thick](3.8,-1)--(3.2,-1.4);
\draw[thick](3.8,-1)--(2.5,0.8);
 \draw[thick](1,-0.5)--(2.5,0.8);
  \draw[thick](3.2,-1.4)--(2.5,0.8);
  \draw[thick](2.5,0.8)--(4,0.3);
    \draw[thick,dashed](2.5,-0.5)--(4,0.3);
        \draw[thick](2.5,-0.5)--(4.2,-0.2);
    \draw[thick](4.2,-0.2)--(4,0.3);
       \draw[thick](4.2,-0.2)--(2.5,0.8);
 \draw[thick,dashed](3.8,-1)--(2.5,-0.5);
 \draw[thick,dashed](1,-0.5)--(2.5,-0.5);
 \node [left ] at (1.5,-1) {$D$};
 \node [left ] at (1,-0.5) {$C$};
 \node [below ] at (2.5,-0.5) {$B$};
 \node [above right ] at (2.5,0.8) {$A$};
  \node [right ] at  (3.2,-1.4) {$E$};
  \node [right ] at (3.8,-1) {$F$};
  \node [right ] at  (4.2,-0.2) {$G$};
  \node [right ] at (4,0.3) {$H$};
  \end{tikzpicture} \qquad
\begin{tikzpicture}
 \node [left] at  (0.7,0) {\emph{ii)}};
 \draw[thick](1,0)--(1.5,-1);
 \draw[thick](1.5,-1)--(2.5,-0.5);
 \draw[thick](2.5,-0.5)--(2.5,0.8);
 \draw[thick](2.5,0)--(2.5,0.8);
 \draw[thick](1.5,-1)--(2.5,0.8);
 \draw[thick](2.5,-0.5)--(3.5,-1);
 \draw[thick](3.7,0)--(3.5,-1);
\draw[thick](3.7,0)--(2.5,0.8);
 \draw[thick](1,0)--(2.5,0.8);
  \draw[thick](3.5,-1)--(2.5,0.8);
 \draw[thick,dashed](3.7,0)--(2.5,-0.5);
 \draw[thick,dashed](1,0)--(2.5,-0.5);
 \node [left ] at (1.5,-1) {$D$};
 \node [above ] at (1,0) {$C$};
 \node [below ] at (2.5,-0.5) {$B$};
 \node [above right ] at (2.5,0.8) {$A$};
  \node [right ] at  (3.5,-1) {$F$};
  \node [right ] at (3.7,0) {$E$};
  \end{tikzpicture} \qquad
  \begin{tikzpicture}
 \node [left] at  (0.7,0) {\emph{iii)}};
 \draw[thick](1,0)--(1.5,-1);
 \draw[thick](1.5,-1)--(2.5,-0.5);
 \draw[thick](2.5,-0.5)--(2.5,0.8);
 \draw[thick](2.5,-0.5)--(2,0.8);
  \draw[thick](2.5,0.8)--(2,0.8);
 \draw[thick](1.5,-1)--(2,0.8);
 \draw[thick](2.5,-0.5)--(3.5,-1);
 \draw[thick](4,0)--(3.5,-1);
\draw[thick](4,0)--(2.5,0.8);
 \draw[thick](1,0)--(2,0.8);
  \draw[thick](3.5,-1)--(2.5,0.8);
 \draw[thick,dashed](4,0)--(2.5,-0.5);
 \draw[thick,dashed](1,0)--(2.5,-0.5);
 \node [left ] at (1.5,-1) {$E$};
 \node [above ] at (1,0) {$D$};
 \node [below ] at (2.5,-0.5) {$B$};
  \node [above left] at (2,0.8) {$C$};
 \node [above right ] at (2.5,0.8) {$A$};
  \node [right ] at  (3.5,-1) {$F$};
  \node [right ] at (4,0) {$G$};
  \end{tikzpicture}
\caption{Dual complexes of the 3-stars in Figure \ref{fig:closedcases}.}\label{fig:closedsubs}
\end{figure}
We claim that the subcomplex are realized as a 3-star on $X$ only if the coefficients $c_{\bf v}$ satisfy some linear conditions. To show this, the idea is to find in each case the equations that the vertex $v$ must satisfy and arrange them in a matrix form, similar to \eqref{eq:MLP}. For example, in case \emph{i)}, we see that $v$ lies on each of the planes spanned by $(AB)^\vee$, $(CD)^\vee$, $(EF)^\vee$ and $(GH)^\vee$. Writing out the corresponding equations, we obtain
\begin{equation}\label{eq:case1}
\begin{pmatrix}
\ora{AB}\;&c_B-c_A\\
\ora{CD}\;&c_D-c_C\\
\ora{EF}\;&c_F-c_E\\
\ora{GH}\;&c_H-c_G
\end{pmatrix}
\begin{pmatrix}
v\\
1
\end{pmatrix}=0.
\end{equation}
Observe that the leftmost matrix in \eqref{eq:case1} is a $4\times 4$-matrix; let us call it $M$. Since the null-space of $M$ is non-trivial (it contains $(v,1)^T$), we must have $\det M=0$, giving a linear relation in the coefficients $c$. Note that this would reduce to an empty condition  if $\rank M\leq 2$, but it is easy to see that in our case $\ora{AB},\ora{CD},\ora{EF},\ora{GH}$ span all of $\rr^3$, so $\rank M=3$. This proves the claim in case \emph{i)}.

The cases \emph{i)} and \emph{ii)} are  done in the same way, but with the matrix $M$ exchanged with 
\begin{equation*}
ii)\;\;M'=\begin{pmatrix}
\ora{AB}\;&c_B-c_A\\
\ora{AC}\;&c_C-c_A\\
\ora{AD}\;&c_D-c_A\\
\ora{EF}\;&c_F-c_E
\end{pmatrix}, \qquad \qquad
iii)\;\;M''=\begin{pmatrix}
\ora{AB}\;&c_B-c_A\\
\ora{AC}\;&c_C-c_A\\
\ora{DE}\;&c_E-c_D\\
\ora{FG}\;&c_G-c_F
\end{pmatrix}.
\end{equation*}

We can now show that a general smooth $X$ contains no  3-stars as above. Indeed, let $\ks$ be any unimodular triangulation of $d\,\Delta_3$, and $c$ a point in the secondary cone $\Sigma(\ks)$. Then as we have seen, any 3-star subcomplex in $\ks$ like those in Figure \ref{fig:closedsubs} can be realized on $X_c$ only if $c$ lies on a certain hyperplane. Moreover, $\ks$ contains at most finitely many of the 3-star subcomplexes in Figure \ref{fig:closedsubs}. Hence for any $c$ in the complement of a finite union of hyperplanes, $X_c$ contains no special 3-stars. 

We conclude by remarking that if the vertices of the surface in case \emph{i)} are adjacent, the tetrahedra in the dual complex might share some faces. This would not influence our arguments above. 
\end{proof}

\begin{remark} The hypothesis that the two vertices $V_1$ and $V_2$ are not adjacent on $X$ is necessary. In fact, suppose that they are adjacent, so there is an edge $E$ of $X$ connecting them. The 3-star is illustrated in Figure \ref{fig:newmotif}, where the egde is $E$ is represented by the dashed line. The dual subcomplex is the one depicted on the right-hand side. The two vertices $V_1$ and $V_2$ are dual  to the tetrahedra $ABCD$ and $ABCF$. They share the face $ABC$ which is dual to the edge $E$. Then the vertex $v$ of the 3-star satisfies the following equations: 
\begin{equation}
\begin{pmatrix}
\ora{BC}\;&c_C-c_B\\
\ora{AD}\;&c_D-c_A\\
\ora{AF}\;&c_F-c_A
\end{pmatrix}
\begin{pmatrix}
v\\
1
\end{pmatrix}=0.
\end{equation}
They are not enough to impose any linear relation among the coefficients. Therefore the surface $X$ is general. When drawing the motif, we  will always add a dotted line connecting the two vertices $V_1$ and $V_2$ in order to emphasize that they are adjacent. 
\begin{figure}[h]
\begin{tikzpicture}
\draw [thick] (-1,1)--(0,0);
\draw [thick] (-1,-1)--(0,0);
\draw [thick] (0,0) --(2,0);
\draw [thick, dashed] (-0.5,0.5)--(-0.5,-0.5);
 \node [above] at (-0.5,0.6) {$V_1$};
  \node [below] at (-0.5,-0.6) {$V_2$};
\draw [fill] (-0.5,0.5) circle [radius=.1];
\draw [fill] (-0.5,-0.5) circle [radius=.1];
\draw [ultra thick] (0,0.5)--(0,-0.5);
 \node [above right] at (0,0) {$v$};
  \node [left] at (-0.5,0) {$E$};
 \end{tikzpicture}
 \,\,\,\,
 \begin{tikzpicture}
 \draw [thick,dashed] (0,0)--(1,0);
\draw [thick,dashed] (-1,1)--(0,0);
\draw [thick](-1,1)--(1,0);
\draw[thick](-1.5,-1)--(1,0);
\draw [thick,dashed] (-1.5,-1) --(0,0);
\draw [thick,dashed] (-2,0.5) --(0,0);
\draw [thick,dashed] (-2,0.5) --(1,0);
\draw [thick] (-1.5,-1)--(-1,1);
\draw [thick](-2,0.5)--(-1,1);
\draw[thick](-1.5,-1)--(-2,0.5);
 \node  [above] at (-1,1) {$A$};
 \node [right] at (1,0) {$B$};
 \node [below] at (0,0) {$C$};
 \node [left] at (-2,0.5) {$D$};
 \node [left] at (-1.5,-1) {$F$};
 \end{tikzpicture}

\caption{A special 3-star and its dual complex.}\label{fig:newmotif}
\end{figure}
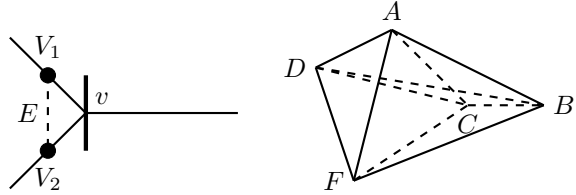
\end{remark}
From Lemma \ref{lem:special} and Lemma \ref{lem:doublytres} it is now possible to construct a list containing all possible motifs of a non-degenerate tropical line on a general smooth tropical surface. The result is shown in Table \ref{bigtable}. 

\begin{table}
\centering
\includegraphics[angle=90,height=0.9\textheight]{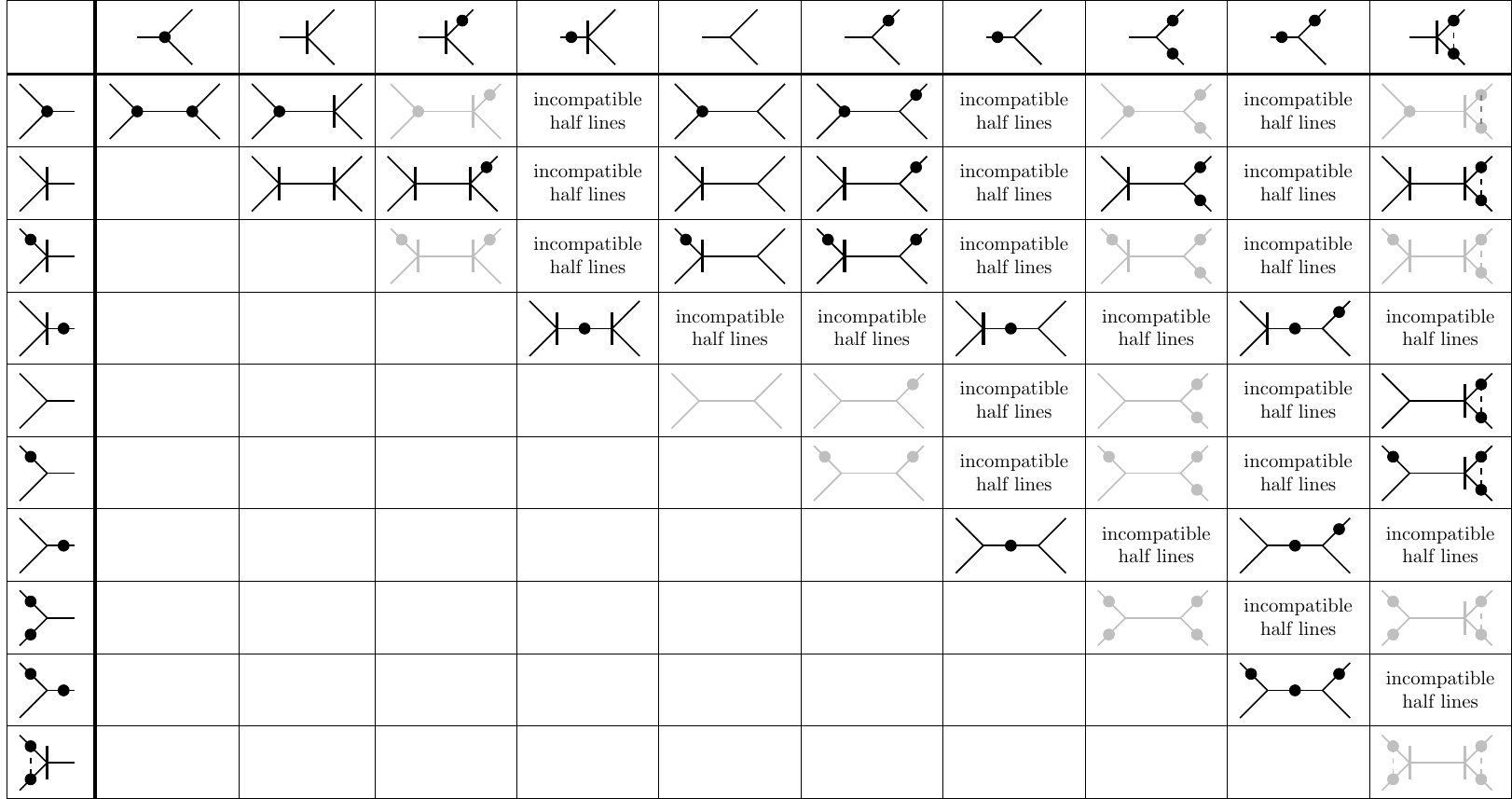}\\

\caption{Motifs of tropical lines on a general smooth surface. Gray=impossible or non-general.}\label{bigtable}
\end{table}

Note that we do not claim that all the entries of Table \ref{bigtable} actually occur on general smooth surfaces. For starters, the following motifs are  impossible on any tropical surface: The middle segment of $L$ is contained in a 2-cell of $X$, but the part of $L$ contained in this cell spans $\rr^3$, which is a contradiction.
\begin{center}
\begin{tikzpicture}
\draw [thick] (-0.5,0.5)--(0,0);
\draw [thick] (-0.5,-0.5)--(0,0);
\draw [thick] (0,0) --(1,0);
\draw [thick] (1,0) --(1.5,0.5);
\draw [thick] (1,0) --(1.5,-0.5);
\end{tikzpicture} 
\,\,\,
\begin{tikzpicture}
\draw [thick] (-0.5,0.5)--(0,0);
\draw [thick] (-0.5,-0.5)--(0,0);
\draw [thick] (0,0) --(1,0);
\draw [thick] (1,0) --(1.5,0.5);
\draw [thick] (1,0) --(1.5,-0.5);
\draw [fill] (1.3,0.3) circle [radius=.07];
\end{tikzpicture}
\,\,\,
\begin{tikzpicture}
\draw [thick] (-0.5,0.5)--(0,0);
\draw [thick] (-0.5,-0.5)--(0,0);
\draw [thick] (0,0) --(1,0);
\draw [thick] (1,0) --(1.5,0.5);
\draw [thick] (1,0) --(1.5,-0.5);
\draw [fill] (1.3,0.3) circle [radius=.07];
\draw [fill] (1.3,-0.3) circle [radius=.07];
\end{tikzpicture}
\,\,\,
\begin{tikzpicture}
\draw [thick] (-0.5,0.5)--(0,0);
\draw [thick] (-0.5,-0.5)--(0,0);
\draw [thick] (0,0) --(1,0);
\draw [thick] (1,0) --(1.5,0.5);
\draw [thick] (1,0) --(1.5,-0.5);
\draw [fill] (1.3,0.3) circle [radius=.07];
\draw [fill] (-0.3,0.3) circle [radius=.07];
\end{tikzpicture}
\,\,\,
\begin{tikzpicture}
\draw [thick] (-0.5,0.5)--(0,0);
\draw [thick] (-0.5,-0.5)--(0,0);
\draw [thick] (0,0) --(1,0);
\draw [thick] (1,0) --(1.5,0.5);
\draw [thick] (1,0) --(1.5,-0.5);
\draw [fill] (1.3,0.3) circle [radius=.07];
\draw [fill] (1.3,-0.3) circle [radius=.07];
\draw [fill] (-0.3,0.3) circle [radius=.07];
\end{tikzpicture}
\,\,\,
\begin{tikzpicture}
\draw [thick] (-0.5,0.5)--(0,0);
\draw [thick] (-0.5,-0.5)--(0,0);
\draw [thick] (0,0) --(1,0);
\draw [thick] (1,0) --(1.5,0.5);
\draw [thick] (1,0) --(1.5,-0.5);
\draw [fill] (1.3,0.3) circle [radius=.07];
\draw [fill] (1.3,-0.3) circle [radius=.07];
\draw [fill] (-0.3,0.3) circle [radius=.07];
\draw [fill] (-0.3,-0.3) circle [radius=.07];
\end{tikzpicture}
\end{center}

Furthermore, we have the following lemma: 
\begin{lemma}\label{lem:non-gen}
\begin{enumerate}[i)]
\item A general smooth $X$ has no tropical lines with motifs shown in Figure \ref{fig:non-gen}.
\item A general smooth $X$ has no tropical line such that its motif is the one in Figure \ref{fig:nongeneric} and its middle segment goes across a 2-dimensional cell of $X$.
\end{enumerate}
\end{lemma}

\begin{figure}[!htb]
    \centering
    \begin{minipage}{.6\textwidth}
        \centering
        \begin{tikzpicture}
\draw [thick] (-0.5,0.5)--(0,0);
\draw [thick] (-0.5,-0.5)--(0,0);
\draw [thick] (0,0) --(1,0);
\draw [thick] (1,0) --(1.5,0.5);
\draw [thick] (1,0) --(1.5,-0.5);
\draw [fill] (0,0) circle [radius=.07];
\draw [fill] (1.3,0.3) circle [radius=.07];
\draw [fill] (1.3,-0.3) circle [radius=.07];
\end{tikzpicture}
\,
        \begin{tikzpicture}
\draw [thick] (-0.5,0.5)--(0,0);
\draw [thick] (-0.5,-0.5)--(0,0);
\draw [thick] (0,0) --(1,0);
\draw [thick] (1,0) --(1.5,0.5);
\draw [thick] (1,0) --(1.5,-0.5);
\draw [fill] (0,0) circle [radius=.07];
\draw [fill] (1.3,0.3) circle [radius=.07];
\draw [ultra thick] (1,0.3)--(1,-0.3);
\end{tikzpicture}
\,
        \begin{tikzpicture}
\draw [thick] (-0.5,0.5)--(0,0);
\draw [thick] (-0.5,-0.5)--(0,0);
\draw [thick] (0,0) --(1,0);
\draw [thick] (1,0) --(1.5,0.5);
\draw [thick] (1,0) --(1.5,-0.5);
\draw [fill] (-0.3,0.3) circle [radius=.07];
\draw [fill] (1.3,0.3) circle [radius=.07];
\draw [ultra thick] (1,0.3)--(1,-0.3);
\draw [ultra thick] (0,0.3)--(0,-0.3);
\end{tikzpicture}
\,
 \begin{tikzpicture}
\draw [thick] (-0.5,0.5)--(0,0);
\draw [thick] (-0.5,-0.5)--(0,0);
\draw [thick] (0,0) --(1,0);
\draw [thick] (1,0) --(1.5,0.5);
\draw [thick] (1,0) --(1.5,-0.5);
\draw [fill] (-0.3,0.3) circle [radius=.07];
\draw [fill] (1.3,-0.3) circle [radius=.07];
\draw [fill] (1.3,0.3) circle [radius=.07];
\draw [ultra thick] (0,0.3)--(0,-0.3);
\end{tikzpicture}

\vspace{\baselineskip}

        \begin{tikzpicture}
\draw [thick] (-0.5,0.5)--(0,0);
\draw [thick] (-0.5,-0.5)--(0,0);
\draw [thick] (0,0) --(1,0);
\draw [thick] (1,0) --(1.5,0.5);
\draw [thick] (1,0) --(1.5,-0.5);
\draw [fill] (0,0) circle [radius=.07];
\draw [fill] (1.3,0.3) circle [radius=.07];
\draw [fill] (1.3,-0.3) circle [radius=.07];
\draw [ultra thick] (1,0.3)--(1,-0.3);
\draw [thick, dashed] (1.3,0.3)--(1.3,-0.3);
\end{tikzpicture}
\,
        \begin{tikzpicture}
\draw [thick] (-0.5,0.5)--(0,0);
\draw [thick] (-0.5,-0.5)--(0,0);
\draw [thick] (0,0) --(1,0);
\draw [thick] (1,0) --(1.5,0.5);
\draw [thick] (1,0) --(1.5,-0.5);
\draw [fill] (-0.3,0.3) circle [radius=.07];
\draw [fill] (1.3,0.3) circle [radius=.07];
\draw [fill] (1.3,-0.3) circle [radius=.07];
\draw [ultra thick] (1,0.3)--(1,-0.3);
\draw [ultra thick] (0,0.3)--(0,-0.3);
\draw [thick, dashed] (1.3,0.3)--(1.3,-0.3);
\end{tikzpicture}
\,
        \begin{tikzpicture}
\draw [thick] (-0.5,0.5)--(0,0);
\draw [thick] (-0.5,-0.5)--(0,0);
\draw [thick] (0,0) --(1,0);
\draw [thick] (1,0) --(1.5,0.5);
\draw [thick] (1,0) --(1.5,-0.5);
\draw [fill] (-0.3,0.3) circle [radius=.07];
\draw [fill] (-0.3,-0.3) circle [radius=.07];
\draw [fill] (1.3,0.3) circle [radius=.07];
\draw [fill] (1.3,-0.3) circle [radius=.07];
\draw [ultra thick] (1,0.3)--(1,-0.3);
\draw [thick, dashed] (1.3,0.3)--(1.3,-0.3);
\end{tikzpicture}
\,
        \begin{tikzpicture}
\draw [thick] (-0.5,0.5)--(0,0);
\draw [thick] (-0.5,-0.5)--(0,0);
\draw [thick] (0,0) --(1,0);
\draw [thick] (1,0) --(1.5,0.5);
\draw [thick] (1,0) --(1.5,-0.5);
\draw [fill] (-0.3,0.3) circle [radius=.07];
\draw [fill] (-0.3,-0.3) circle [radius=.07];
\draw [fill] (1.3,0.3) circle [radius=.07];
\draw [fill] (1.3,-0.3) circle [radius=.07];
\draw [ultra thick] (1,0.3)--(1,-0.3);
\draw [ultra thick] (0,0.3)--(0,-0.3);
\draw [thick, dashed] (1.3,0.3)--(1.3,-0.3);
\draw [thick, dashed] (-0.3,0.3)--(-0.3,-0.3);
\end{tikzpicture}
            \caption{Motifs of Lemma \ref{lem:non-gen}\,\emph{i)}.}\label{fig:non-gen}
    \end{minipage}%
    \begin{minipage}{0.4\textwidth}
        \centering
          \begin{tikzpicture}
\draw [thick] (-0.5,0.5)--(0,0);
\draw [thick] (-0.5,-0.5)--(0,0);
\draw [thick] (0,0) --(1,0);
\draw [thick] (1,0) --(1.5,0.5);
\draw [thick] (1,0) --(1.5,-0.5);
\draw [fill] (0,0) circle [radius=.07];
\draw [fill] (1,0) circle [radius=.07];
\end{tikzpicture}
\caption{Motif of Lemma \ref{lem:non-gen}\,\emph{ii)}.}
 \label{fig:nongeneric}
    \end{minipage}
    \end{figure}

%

\begin{proof}
\emph{i)} The idea is basically the same as in the proofs of Lemma \ref{lem:doublytres} and Lemma \ref{lem:special}. Each case implies some linear relations between the coefficients $c_{\bf v}$ of the polynomial defining $X$. We sketch a typical argument: Suppose $L$ has the top left motif, and observe that the line subcomplex of $L$ is then homeomorphic to case \emph{i)} of Figure \ref{fig:closedsubs}. Let $v_1$ be the vertex of $L$ which is also a vertex of $X$; assume this is dual to the tetrahedron $ABGH$. Then, by the property of duality, $v_1$ is uniquely determined by $c_A, c_B, c_G$ and $c_H$. (In fact, the coordinates of $v_1$ are linear forms in these $c$'s.) Similarly, the other vertex of $L$, $v_2$, is determined by $c_A, c_B, c_C, c_D, c_E$ and $c_F$ (this corresponds to solving the equation \eqref{eq:case1}, but with the last row removed.) Finally, since $v_1v_2$ is the middle segment of a tropical line, it has a prescribed direction. This forces a linear relation between the coefficients.
The remaining  cases of the above claim are done similarly. For each motif we can see that each vertex of the line need to satisfies three linear equations. Moreover the middle segment has a prescribed direction, so there is a linear relation between the coefficients. 

\emph{ii)} This motif was discussed in Remark \ref{rem:middleedge}. If the middle segment of $L$ goes across a $2$-dimensional cell of $X$, the line subcomplex is homeomorphic to case \emph{ii)} of Figure \ref{fig:tvecp}. In this case the argument sketched in \emph{i)} applies again: Each vertex of $L$ is determined by the coefficients, and the direction vector of the middle segment implies a linear relation between these.

Note that this argument does not apply if the middle segment of $L$ is a $1$-dimensional cell of $X$: In this case the direction of the middle segment is encoded in the line subcomplex as a normal vector of the common facet of the two tetrahedra as in case \emph{i)} of Figure \ref{fig:tvecp}.
\end{proof}

\section{The classification theorem}
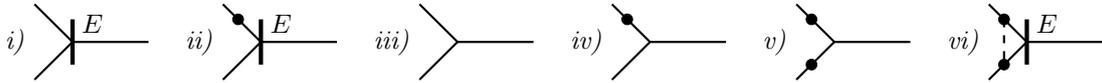
\begin{figure}[h]
\begin{tikzpicture}
\node [left] at (-0.5,0) {\emph{i)}};
\draw [thick] (-0.5,0.5)--(0,0);
\draw [thick] (-0.5,-0.5)--(0,0);
\draw [thick] (0,0) --(1,0);
\draw [ultra thick] (0,0.3)--(0,-0.3);
 \node [above right] at (0,0) {$E$};
 \end{tikzpicture}
 \ \
   \begin{tikzpicture}
   \node [left] at (-0.5,0) {\emph{ii)}};
\draw [thick] (-0.5,0.5)--(0,0);
\draw [thick] (-0.5,-0.5)--(0,0);
\draw [thick] (0,0) --(1,0);
\draw [fill] (-0.3,0.3) circle [radius=.07];
\draw [ultra thick] (0,0.3)--(0,-0.3);
 \node [above right] at (0,0) {$E$};
 \end{tikzpicture}
 \ \
 \begin{tikzpicture}
 \node [left] at (-0.5,0) {\emph{iii)}};
\draw [thick] (-0.5,0.5)--(0,0);
\draw [thick] (-0.5,-0.5)--(0,0);
\draw [thick] (0,0) --(1,0);
 \end{tikzpicture}
 \ \
  \begin{tikzpicture}
  \node [left] at (-0.5,0) {\emph{iv)}};
\draw [thick] (-0.5,0.5)--(0,0);
\draw [thick] (-0.5,-0.5)--(0,0);
\draw [thick] (0,0) --(1,0);
\draw [fill] (-0.3,0.3) circle [radius=.07];
 \end{tikzpicture}
 \ \
  \begin{tikzpicture}
  \node [left] at (-0.5,0) {\emph{v)}};
\draw [thick] (-0.5,0.5)--(0,0);
\draw [thick] (-0.5,-0.5)--(0,0);
\draw [thick] (0,0) --(1,0);
\draw [fill] (-0.3,0.3) circle [radius=.07];
\draw [fill] (-0.3,-0.3) circle [radius=.07];
 \end{tikzpicture}
 \ \
 \begin{tikzpicture}
 \node [left] at (-0.5,0) {\emph{vi)}};
\draw [thick] (-0.5,0.5)--(0,0);
\draw [thick] (-0.5,-0.5)--(0,0);
\draw [thick] (0,0) --(1,0);
\draw [thick, dashed] (-0.3,0.3)--(-0.3,-0.3);
\draw [fill] (-0.3,0.3) circle [radius=.07];
\draw [fill] (-0.3,-0.3) circle [radius=.07];
\draw [ultra thick] (0,0.3)--(0,-0.3);
 \node [above right] at (0,0) {$E$};
 \end{tikzpicture}
\caption{Tropical half lines.}\label{fig:halflines} 
\end{figure}

In the last section we identified 14 entries of Table \ref{bigtable} that were either impossible on any tropical surface $X$, or {\em non-general}, meaning that they do not occur on general smooth $X$. In this section we analyze the remaining 20 motifs. To avoid repeating ourselves too much, we start by giving some auxiliary observations about {\em tropical half lines} on $X$, which we will apply frequently. A tropical half line in $\rr^3$ is the remaining part of a non-degenerate tropical line, after removing two adjacent rays. Figure \ref{fig:halflines} shows tropical half lines on $X$ in different positions.

For a tropical half line $H$, let $H^b$ be its bounded segment. Note that if $H^b\sub X$ is non-trespassing, then there is a unique cell of $X$, denoted $C^b$, containing $H^b$. The following lemma gives information on the position of the dual cell $(C^b)^\vee$ in the triangulation $ \ks_X$. As always, $X$ is assumed to be smooth of some fixed degree $d$. 
 
\begin{lemma}\label{halfb}
Let $H\sub X$ be a tropical half line with unbounded rays in the directions $\omega_i$ and $\omega_j$, such that $H^b$ is non-trespassing and contained in a cell $C^b$ of $X$. Then $(C^b)^\vee$ is contained in a plane with equation $x_i+x_j=K$ for some non-negative integer $K$.
\end{lemma}
\begin{proof}
Recall that any vector contained in $C^b$ is orthogonal to $(C^b)^\vee$. In the case where $i,j\neq 0$, this immediately proves the assertion, since by assumption $C^b$ contains the vector $-(\omega_i+\omega_j)=-(e_i+e_j)$. For the remaining case, suppose $j=0$, and let $(i,i',i'')$ be any permutation of $(1,2,3)$. Then $C^b$ contains the vector $\omega_i+\omega_0=-(e_{i'}+e_{i''})$, so $(C^b)^\vee$ lies in a plane with equation $x_{i'}+x_{i''}=\text{constant}$. This is equivalent to the statement in the lemma, since $x_i+x_0=K\Longleftrightarrow  x_{i'}+x_{i''}=d-K$.
\end{proof}
\begin{lemma}\label{halflines}
Let $H$ be as in Lemma \ref{halfb}. Let $v$ be the vertex of $H$, and $k$ is the number of unbounded rays of $H$ which are trespassing on $X$. If either
\begin{itemize}
\item $\dim(\emph{\mc}(v))\geq k$ and $\dim(\emph{\mc}(v))>0$, or 
\item $H$ is the half line vi) in Figure \ref{halfb}, or 
\item $\dim(\emph{\mc}(v))\geq k$,  $\dim(\emph{\mc}(v))=0$ and $\dim(C^b)=1$,

\end{itemize}
then $(C^b)^\vee$ lies in the plane with equation $x_i+x_j=\dim(C^b)-\dim(\emph{\mc}(v))+k$. 
\end{lemma}

\begin{proof}
By symmetry we can assume that $i=1$ and $j=2$, so by Lemma~\ref{halfb}, $(C^b)^\vee$ lies in a plane given by $x_1+x_2=K$ for some $K$. For  $\dim(\mc(v))>0$, it is easy to see that it suffices to consider the five first cases shown in Figure~\ref{fig:halflines}. Note that in all these cases, $\dim( C^b)=2$.
\smallskip 

In case \emph{i)}, we have $\dim(\mc(v))=1$ and $k=0$, so we must show that $K=2+0-1=1$. We see that the triangle $E^\vee$ has one edge on $F_1$, another edge on $F_2$, while its last edge is $(C^b)^\vee$. Hence the vertices of $E^\vee$ are of the form $(0,0,a)$, $(0,K,b)$ and $(K,0,c)$, where $a,b,c\in\nn_0$. Let $P=(p,q,r)$ be the fourth vertex of a tetrahedron in $\ks_X$ having $E^\vee$ as a facet. A standard calculation shows that the volume of this tetrahedron is divisible by $K$. Hence unimodularity of $\ks_X$ implies $K=1$, as wanted.
\smallskip 

In case \emph{ii)}, we must show that $K=2$. Here, $E^\vee$ has one edge in $F_2$, another in the plane $x_1=1$, and the third is $(C^b)^\vee$. Thus the vertices of $E^\vee$ are $(1,0,a),(K,0,b),(1,K-1,c)$ for some $a,b,c\in\nn_0$. A volume calculation shows that any unimodular tetrahedron having $E^\vee$ as facet, has a volume divisible by $K-1$. Thus $K=2$.
\smallskip

In case \emph{iii)}, we must show $K=0$. It is clear that $(C^b)^\vee$ lies in both facets $F_1$ and $F_2$, and therefore in the edge $F_{12}$ of $d \, \Delta_3$. Since $F_{12}$ is contained in the plane $x_1+x_2=0$, we are done.
\smallskip

In case \emph{iv)}, we find similarly that $(C^b)^\vee$ lies in the intersection of $F_2$ (where $x_2=0$) and the plane given by $x_1=1$. In particular, $(C^b)^\vee$ lies in the plane where $x_1+x_2=1$, as claimed in the lemma.
\smallskip

In case \emph{v)}, the edge $(C^b)^\vee$ lies in the intersection of the planes with equations $x_1=1$ and $x_2=1$. In particular, this means $x_1+x_2=2$, which is again what we needed to prove.

\smallskip

Finally, in case \emph{vi)}, the triangle $E^\vee$ has one edge on $x_1 =1$, one edge on $x_2 = 1$ and the third edge on $x_1 + x_2 = K$. Therefore its vertices are $A = (1, K-1, a)$, $B= (1,1,b)$ and $C=(K-1,1,c)$. Again a volume computation shows that the volume of a tetrahedron having $E^\vee$ as face is divisible by $K-2$. Unimodularity implies that $K =3$. 

\smallskip

It remains to treat the case where $\dim(\mc(v))=0$ and $\dim(C^b)=1$. In other words, $v$ is a vertex of $X$ and $H^b$ is contained in an edge $C^b$ of $X$. Dually, $(C^b)^\vee$ is a triangle $ABC\in \ks_X$, and $v^\vee$ is a tetrahedron having $(C^b)^\vee$ as a facet, i.e. $v^\vee=ABCD$ for some lattice  point $D$. By Lemma~\ref{halfb}, $ABC$ lies in a plane with equation $x_1+x_2=K$, and by  volume considerations, $D$ must then lie in a plane given by $x_1+x_2=K\pm1$. We have to prove that $K=1$. 

Observe that $ABCD$ has exits in both directions $\omega_1$ and $\omega_2$, since $H$ has no trespassing rays given  the assumption $k\leq \dim(\mc(v))$. Now, if the triangle $ABC$ has exits in neither of the two directions, then we must have $D\in F_{12}$, implying $K=1$. If $ABC$ has an exit in exactly one of the two directions, say $\omega_1$, then we must have (possibly after renaming) that $AB\sub F_1$ and $CD\sub F_2$. Writing out what this means for the coordinates of $A$, $B$, $C$ and $D$, a volume calculation of $ABCD$ shows that unimodularity again implies $K=1$. Note that $ABC$ (being a non-degenerate triangle contained in the plane of the form $x_1+x_2=K$) cannot have exits in both directions $\omega_1$ and $\omega_2$. Hence we have covered all cases. 
\end{proof}
\begin{remark}
Lemma~\ref{halflines} implies in particular that in the cases mentioned, the integer $K$ in Lemma~\ref{halfb} is dependent only on $\dim(C^b)$, $\dim(\mc(v))$ and $k$. As we will see, this does no longer hold if $\dim(\mc(v))=0$ and $\dim(C^b)=2$ (the only case not covered by the lemma), and this fact is what allows the motifs 3G and 3H to occur on surfaces of arbitrarily high degree.
\end{remark}


\begin{table}
\centering
\begin{tabular}{|c|c|c|c|} 
 \hlineB{3.5}
\multicolumn{4}{!{\vrule width 2pt}c!{\vrule width 2pt}}{Primal motifs of isolated lines}\\
 \hlineB{3.5}
\raisebox{0.5cm}{\multirow{1}{*}{$\textrm{deg}(X) = 1$}}& \begin{tikzpicture} \smallskeleton  \node [above] at (0.5,0.5){1A}; \draw [ultra thick] (0,0.3)--(0,-0.3); \end{tikzpicture} &   \begin{tikzpicture} \smallskeleton  \node [above] at (0.5,0.5){1B}; \draw [fill] (0.5,0) circle [radius=.1];  \end{tikzpicture}  & \\
\hline
 \multirow{2}{*}{$\textrm{deg}(X) = 2$} & \begin{tikzpicture} \smallskeleton  \node [above] at (0.5,0.5){2A}; \draw [ultra thick] (0,0.3)--(0,-0.3); \draw [fill] (-0.3,0.3) circle [radius=.1];  \end{tikzpicture} &  \begin{tikzpicture} \smallskeleton  \node [above] at (0.5,0.5){2B}; \draw [ultra thick] (0,0.3)--(0,-0.3); \draw [fill] (1.3,0.3) circle [radius=.1];  \end{tikzpicture} & \begin{tikzpicture} \smallskeleton  \node [above] at (0.5,0.5){2C};  \draw [ultra thick] (0,0.3)--(0,-0.3); \draw [fill] (0.5,0) circle [radius=.1]; \end{tikzpicture}\\ \cline{2-4}
 & \begin{tikzpicture} \smallskeleton  \node [above] at (0.5,0.5){2D};  \draw [ultra thick] (0,0.3)--(0,-0.3); \draw [ultra thick] (1,0.3)--(1,-0.3); \end{tikzpicture} &  \begin{tikzpicture} \smallskeleton  \node [above] at (0.5,0.5){2E};   \draw [fill] (0.5,0) circle [radius=.1];
\draw [fill] (1.3,0.3) circle [radius=.1]; \end{tikzpicture}  & \begin{tikzpicture} \smallskeleton  \node [above] at (0.5,0.5){2F};  \draw [fill] (0,0) circle [radius=.1]; \draw [fill] (1,0) circle [radius=.1]; \end{tikzpicture} \\
\hline
\raisebox{0cm}{\multirow{3}{*}{$\textrm{deg}(X) = 3$}}  & \begin{tikzpicture} \smallskeleton  \node [above] at (0.5,0.5){3A};\draw [ultra thick] (0,0.3)--(0,-0.3); \draw [ultra thick] (1,0.3)--(1,-0.3); \draw [fill] (1.3,0.3) circle [radius=.1]; \end{tikzpicture} & \begin{tikzpicture} \smallskeleton  \node [above] at (0.5,0.5){3B}; \draw [ultra thick] (0,0.3)--(0,-0.3); \draw [ultra thick] (1,0.3)--(1,-0.3); \draw [fill] (0.5,0) circle [radius=.1];\end{tikzpicture} & \begin{tikzpicture} \smallskeleton  \node [above] at (0.5,0.5){3C}; \draw [ultra thick] (0,0.3)--(0,-0.3);
\draw [fill] (0.5,0) circle [radius=.1]; \draw [fill] (1.3,0.3) circle [radius=.1];\end{tikzpicture}  \\  \cline{2-4}
 & \begin{tikzpicture} \smallskeleton  \node [above] at (0.5,0.5){3D};\draw [ultra thick] (0,0.3)--(0,-0.3);
\draw [fill] (-0.3,0.3) circle [radius=.1]; \draw [fill] (1.3,0.3) circle [radius=.1]; \end{tikzpicture}  & \begin{tikzpicture} \smallskeleton  \node [above] at (0.5,0.5){3E}; \draw [ultra thick] (0,0.3)--(0,-0.3); \draw [fill] (1.3,-0.3) circle [radius=.1]; \draw [fill] (1.3,0.3) circle [radius=.1]; \end{tikzpicture}  &\begin{tikzpicture} \smallskeleton  \node [above] at (0.5,0.5){3F};  \draw [fill] (0.5,0) circle [radius=.1];
\draw [fill] (1.3,0.3) circle [radius=.1]; \draw [fill] (-0.3,0.3) circle [radius=.1];\end{tikzpicture} \\  \cline{2-4}
\hline
\raisebox{0.5cm}{$\textrm{deg}(X) \geq 2$} & \begin{tikzpicture} \smallskeleton  \node [above] at (0.5,0.5){3G}; \draw [fill] (0,0) circle [radius=.1]; \draw [fill] (1.3,0.3) circle [radius=.1]; \end{tikzpicture} &  \begin{tikzpicture} \smallskeleton  \node [above] at (0.5,0.5){3H}; \draw [fill] (0,0) circle [radius=.1]; \draw [ultra thick] (1,0.3)--(1,-0.3); \end{tikzpicture}& \\
\hline
\raisebox{0.5cm}{$\textrm{deg}(X) = 4$} &\begin{tikzpicture} \smallskeleton \node [above] at (0.5,0.5){4A};\draw [fill] (-0.3,0.3) circle [radius=.1]; \draw [fill] (-0.3,-0.3) circle [radius=.1]; \draw [ultra thick] (0,0.3)--(-0,-0.3); \draw [thick, dashed] (-0.3,0.3)--(-0.3,-0.3); \draw [fill] (1.3,0.3) circle [radius=.1];\end{tikzpicture}  &  \begin{tikzpicture} \smallskeleton  \node [above] at (0.5,0.5){4B};\draw [fill] (-0.3,0.3) circle [radius=.1]; \draw [fill] (-0.3,-0.3) circle [radius=.1]; \draw [ultra thick] (0,0.3)--(-0,-0.3); \draw [thick, dashed] (-0.3,0.3)--(-0.3,-0.3); \draw [ultra thick] (1,0.3)--(1,-0.3);\end{tikzpicture}  & \\
\hline 
\multicolumn{4}{c}{}\\
 \hlineB{3.5}
\multicolumn{4}{!{\vrule width 2pt}c!{\vrule width 2pt}}{Primal motifs of families of lines}\\
 \hlineB{3.5}
\raisebox{0.5cm}{$\textrm{deg}(X) \geq 1$} &\begin{tikzpicture} \smallskeleton  \node [above] at (0.5,0.5){3I}; \draw [fill] (0,0) circle [radius=.1];\end{tikzpicture}  &   & \\
\hline
 \raisebox{0.5cm}{$\textrm{deg}(X) =3$} &\begin{tikzpicture} \smallskeleton  \node [above] at (0.5,0.5){3J};\draw [fill] (-0.3,0.3) circle [radius=.1];
\draw [fill] (-0.3,-0.3) circle [radius=.1]; \draw [ultra thick] (0,0.3)--(-0,-0.3); \draw [thick, dashed] (-0.3,0.3)--(-0.3,-0.3);\end{tikzpicture}  &  & \\ \cline{2-4}
\hline
\end{tabular}
\vspace{\baselineskip}
\caption{Classification of general primal motifs for isolated lines and two-point families.} \label{tab:classif}
\end{table}

\begin{proposition} \label{prop:new} Let $X$ be a general smooth tropical  surface.
\begin{enumerate}[a)]
\item The motifs \emph{1A} and \emph{1B} occur for isolated lines  on surfaces of degree $1$.
\item The motifs \emph{2A}, \emph{2B}, \emph{2C}, \emph{2D}, \emph{2E} and \emph{2F}  occur for isolated lines on surfaces of degree $2$. 
\item The motifs \emph{3A}, \emph{3B}, \emph{3C}, \emph{3D}, \emph{3E}, and \emph{3F} occur for isolated lines on surfaces of degree $3$. 
\item The motifs \emph{3G} and \emph{3H} occur for isolate lines on surfaces of any degree $d\geq 2$.
\item The motifs \emph{4A} and \emph{4B} occur  for isolated lines on surfaces of degree $4$.
\item The motif \emph{3I} occurs for families of lines on surfaces of any degree $d\geq 1$.
\item The motif \emph{3J} occurs for families of lines on surfaces of degree $3$. 

\end{enumerate}

\end{proposition}
We call these  20 motifs {\em general motifs}. The motifs 1A, 1B, $\dots$, 3I were presented in the classification carried out by the first author in \cite{Vigeland08}. The motif 3J was founded by Hampe and Joswig in \cite{HampeJoswig17}. We decided to keep the labeling introduced in \cite{Vigeland08} in order to emphasize the historical development of the project. 

\begin{proof} Motif 1A: Suppose $L\sub X$, where $X$ has degree $d$, and $L$ has motif 1A on $X$. We can assume w.l.o.g. that $L$ has combinatorial type $((12)(30))$, so the situation is as follows:

\begin{center}
\begin{tikzpicture}
\node [above] at (0.5,0.5){1A};
\smallskeleton
\draw [ultra thick] (0,0.3)--(0,-0.3);
\node [above left] at (-0.5,0.5) {$1$};
\node [below left] at (-0.5,-0.5) {$2$};
\node [above right] at (1.5,0.5) {$3$};
\node [below right] at (1.5,-0.5) {$0$};
\end{tikzpicture}
\end{center}

Regard $L$ as the union of two tropical half lines on $X$, sharing the same bounded segment. Let $C^b$ be the 2-dimensional cell of $X$ containing this bounded segment. Applying Lemma \ref{halflines} to the half line  with rays 1 and 2, it follows that $(C^b)^\vee$ lies in the plane with equation $x_1+x_2=1$. On the other hand, the same lemma applied to the half line with rays 3 and 0 implies that $(C^b)^\vee$ lies in the plane with equation $x_3+x_0=0$, i.e. $x_1+x_2=d$. We conclude that $d=1$. 

\smallskip

\noindent An analogous argument works in all cases mentioned in  parts  \emph{a), b), c)} and \emph{d)} where the middle segment of $L$ is not trespassing on $X$, i.e. motifs 2A, 2B, 2D, 2F, 3A, 3D, 3E, 3J, 4A and 4B. Note in particular that in case 2F, we have $\dim(C^b)=1$ by Lemma \ref{lem:non-gen}\,\emph{b)}, so Lemma \ref{halflines} applies.

\smallskip

\noindent  Motif 1B: Again we assume that $L$ has combinatorial type $((12)(30))$. 

\begin{center}
\begin{tikzpicture}
\node [above] at (0.5,0.5){1B};
\draw [thick] (-0.5,0.5)--(0,0);
\draw [thick] (-0.5,-0.5)--(0,0);
\draw [thick] (0,0) --(1,0);
\draw [thick] (1,0) --(1.5,0.5);
\draw [thick] (1,0) --(1.5,-0.5);
\draw [fill] (0.5,0) circle [radius=.1];
\node [below] at (0.5,-0.1) {$V$};
\node [above left] at (-0.5,0.5) {$1$};
\node [below left] at (-0.5,-0.5) {$2$};
\node [above right] at (1.5,0.5) {$3$};
\node [below right] at (1.5,-0.5) {$0$};
\end{tikzpicture}
\end{center}

This time we regard $L$ as the union of two tropical half lines on $X$, intersecting in the point $V$. Let $C^b_1$ be the 2-dimensional cell containing the bounded segment of the half line with rays 1 and 2, and similarly $C^b_2$ the 2-dimensional cell containing the bounded segment of the  half line with rays 3 and 0.

Now we apply Lemma \ref{halflines} twice, to find that $(C^b_1)^\vee$ and $(C^b_2)^\vee$ lie in the planes with equations $x_1+x_2=0$, and $x_1+x_2=d$ respectively. But $(C^b_1)^\vee$ and $(C^b_2)^\vee$ are opposite edges of the unimodular tetrahedron $V^\vee$. This is only possibly if $d=1$. 

\smallskip

\noindent  An analogous argument works in all cases mentioned in parts \emph{a)}, \emph{b)} and \emph{c)} where the middle segment is trespassing, i.e. the motifs 2C, 2E, 3B, 3C and 3F. 

\smallskip

\noindent Motif  3G: Assume $d\geq 2$ (this is necessary for $L\sub X$ to have motif 3G), and consider the complex $\kr_G$ with maximal elements $\tau_1=ABCD$ and $\tau_2=CDEF$, where
\[
\begin{array}{c}
A=(0,0,1),\quad B=(0,1,1),\quad C=(d-1,0,0), \\
D=(d-2,1,0),\quad E=(d-1,1,0),\quad F=(d-1,0,1). \\
\end{array}
\]
It is clear that $\kr_G$ is of type 3G. We claim that if $\ks_X$ contains $\kr_G$, then $X$ contains a line with motif 3G. Indeed, this can be checked directly by examining the shape of the 2-dimensional cell of $X$ dual to the edge $CD$. Furthermore, by applying the techniques described in \cite[Section 4]{Vigeland10}, one can construct regular triangulations of $d\,\Delta_3$ containing $\kr_G$, for all $d\geq 2$.
This proves the assertion in Theorem \ref{theo:classif}\,\emph{d)} concerning position 3G.

\smallskip

\noindent Motif 3H: Let $\kr_H$ be the complex with maximal elements the tetrahedron $ABCD$ and the triangle $CDE$, where
\begin{equation*}
A=(0,0,1),\quad B=(0,1,2),\quad C=(d-1,0,0),\quad D=(d-2,1,1), \quad E=(d-1,1,0).
\end{equation*}
Then $\kr_H$ is of type 3H. Suppose $\ks$ is a regular triangulation of $d\,\Delta_3$ containing $\kr_H$. By examining the shape of the 2-dimensional cell $(CD)^\vee$, one can see that $\kr_H$ is a motif in $\ks$, realizable on $X_c$ for all $c$ in some open, full-dimensional cone in $\Sigma(\ks)$. Finally, as above, one can construct regular triangulations of $d\,\Delta_3$ containing $\kr_H$, for all $d\geq 2$.

\smallskip

\noindent  Motif 3I: Consider the tetrahedron $\tau$ with vertices $(0,0,0)$, $(0,0,1)$, $(1,0,d-1)$ and $(d-1,1,0)$. Clearly, $\tau$ is a complex of type 3I. In the proof of \cite[Theorem 9.2]{Vigeland10} we showed that if $\tau\in \ks_X$, then $X$ has a tropical line with motif  3I. Moreover, we showed that for all $d\in\nn$ there exists an unimodular triangulation of $d\,\Delta_3$ which contains $\tau$.
\end{proof}

After this preparatory work, we can prove  Theorem \ref{theo:classif}. 

\begin{proof}[Proof of Theorem~\ref{theo:classif}]
This follows from the arguments in the proof of Proposition \ref{prop:new}. For the motifs of two-point families we also use convexity of the cells of $X$: Suppose $L\sub X$ has motif 3I or 3J, with vertices $v_1$ and $v_2$, where $C:=\mc_X(v_2)$ has dimension 2. If $\vec{u}$ is the vector from $v_1$ to $v_2$, then convexity of $C$ implies that the tropical line $L_t$ with vertices $v_1$ and $v_t:=v_1+t\vec{u}$ lies on $X$ for all $t\geq 0$. 
\end{proof}

From the discussion in Lemma \ref{halflines} and Proposition \ref{prop:new}, it follows that in order for a motif  to be realized by a line, the conditions imposed by the exits on the dual motif are not enough. In Table \ref{tab:gen_comblines} we have summarized the cell structures of the motifs  3A,$\dotsc$, 3J, including the additional information provided by Lemma \ref{halflines} on the dual motif. We call a  subcomplex in a regular unimodular triangulation of $3\,\Delta_3$ having the cells structure of a motif and satisfying the condition in the rightmost column in Table \ref{tab:gen_comblines} an \emph{occurrence of a motif}. It should be noted that the conditions given in Table \ref{tab:gen_comblines} are only necessary conditions to find occurrences of motifs: If $\ks$ is an regular triangulation of $d\, \Delta_3$, then a subcomplex of $\ks$ of a given type is not necessarily a motif in $\ks$.
\begin{table}[htbp]
\begin{center}
\begin{tabular}{|c|c|l|}
\hlineB{3.5}
\multicolumn{3}{!{\vrule width 2pt}c!{\vrule width 2pt}}{Motifs of isolated lines}\\
\hlineB{3.5}
Primal Motif & Dual Motif & Necessary conditions \\
\hline
\raisebox{-0.1cm}{
\begin{tikzpicture}[scale=0.8] \decskeleton \node [above] at (0.5,0.5){3A}; \draw [ultra thick] (0,0.3)--(0,-0.3);
\draw [ultra thick] (1,0.3)--(1,-0.3);
\draw [fill] (1.3,0.3) circle [radius=.1]; \end{tikzpicture}  }&\begin{tikzpicture}[scale=0.7] \dmA \end{tikzpicture} &   \raisebox{0.7cm}{\begin{tabular}{l}Exits: $AB\subseteq F_i, \ \ \ BD\subseteq F_j,$ \\ $\qquad \ \ \ AC \subseteq F_k, \ \ \ EF \subseteq F_l,$\\ \\$AD \subseteq \mathcal{P}_{x_i+x_j=1},$ \\ 
$CD \subseteq \mathcal{P}_{x_l=1}.$ \\ \end{tabular}} \\
\hline
\begin{tikzpicture}[scale=0.8] \decskeleton \node [above] at (0.5,0.5){3B}; \draw [ultra thick] (0,0.3)--(0,-0.3);
\draw [ultra thick] (1,0.3)--(1,-0.3);
\draw [fill] (0.5,0) circle [radius=.1]; \end{tikzpicture}   &\begin{tikzpicture}[scale=0.7] \dmB \end{tikzpicture}  & \raisebox{0.8cm}{\begin{tabular}{l}Exits: $AB\subseteq F_i, \ \ \ AC\subseteq F_j,$ \\ $ \qquad \ \ \ DF \subseteq F_k, \ \ \ EF \subseteq F_l,$\\ \\ $BC \subseteq \mathcal{P}_{x_i+x_j=1},$\\ $ DE \subseteq \mathcal{P}_{x_k + x_l=1}.$ \\ \end{tabular}}\\
\hline
\begin{tikzpicture}[scale=0.8]  \decskeleton \node [above] at (0.5,0.5){3C}; \draw [ultra thick] (0,0.3)--(0,-0.3);
\draw [fill] (0.5,0) circle [radius=.1];
\draw [fill] (1.3,0.3) circle [radius=.1]; \end{tikzpicture}   &\begin{tikzpicture}[scale=0.7] \dmC \end{tikzpicture}  & \raisebox{0.9cm}{\begin{tabular}{l}Exits: $AB\subseteq F_i, \ \ \ AC\subseteq F_j,$ \\ $ \qquad \ \ \ DE \subseteq F_k, \ \ \ FG \subseteq F_l,$\\ \\$BC \subseteq \mathcal{P}_{x_i+x_j=1},$ \\
$ DE \subseteq \mathcal{P}_{x_l=1}\cap F_k.$ \\ \end{tabular}}\\
\hline
\raisebox{-0.1cm}{\begin{tikzpicture}[scale=0.8] \decskeleton \node [above] at (0.5,0.5){3D}; \draw [ultra thick] (0,0.3)--(0,-0.3);
\draw [fill] (-0.3,0.3) circle [radius=.1];
\draw [fill] (1.3,0.3) circle [radius=.1]; \end{tikzpicture}  }&\begin{tikzpicture}[scale=0.7] \dmD \end{tikzpicture} &\raisebox{0.8cm}{\begin{tabular}{l}Exits: $CE\subseteq F_i, \ \ \ AB\subseteq F_j,$ \\ $ \qquad \ \ \ DE \subseteq F_k, \ \ \ FG \subseteq F_l,$\\ \\$CD \subseteq \mathcal{P}_{x_j=1},$ \\
$DE \subseteq \mathcal{P}_{x_l=1}\cap F_k.$ \\ \end{tabular} }\\
\hline
\raisebox{0.1cm}{\begin{tikzpicture}[scale=0.8] \decskeleton \node [above] at (0.5,0.5){3E}; \draw [ultra thick] (0,0.3)--(0,-0.3);
\draw [fill] (1.3,-0.3) circle [radius=.1];
\draw [fill] (1.3,0.3) circle [radius=.1]; \end{tikzpicture} } &\begin{tikzpicture}[scale=0.7] \dmE \end{tikzpicture} &\raisebox{1cm}{\begin{tabular}{l}Exits: $AB\subseteq F_i, \ \ \ AC\subseteq F_j,$ \\ $ \qquad \ \ \ DE \subseteq F_k, \ \ \ FG \subseteq F_l,$\\ \\$BC \subseteq \mathcal{P}_{x_k=1} \cap \mathcal{P}_{x_l=1}$.  \\ \end{tabular}} \\
\hline
\raisebox{-0.1cm}{\begin{tikzpicture}[scale=0.8] \decskeleton \node [above] at (0.5,0.5){3F}; \draw [fill] (0.5,0) circle [radius=.1];
\draw [fill] (1.3,0.3) circle [radius=.1];
\draw [fill] (-0.3,0.3) circle [radius=.1]; \end{tikzpicture} } &\begin{tikzpicture}[scale=0.7] \dmF \end{tikzpicture} & \raisebox{0.7cm}{\begin{tabular}{l}Exits: $CD\subseteq F_i, \ \ \ AB\subseteq F_j,$ \\ $ \qquad \ \ \ EF \subseteq F_k, \ \ \ GH \subseteq F_l,$\\ \\$CD \subseteq \mathcal{P}_{x_j=1} \cap F_i,$ \\
$ EF \subseteq \mathcal{P}_{x_l=1} \cap F_k$.  \\ \end{tabular}}\\
\hline
\raisebox{-0.1cm}{\begin{tikzpicture}[scale=0.8] \decskeleton \node [above] at (0.5,0.5){3G}; \draw [fill] (0,0) circle [radius=.1];
\draw [fill] (1.3,0.3) circle [radius=.1]; \end{tikzpicture} }&\begin{tikzpicture}[scale=0.7] \dmG \end{tikzpicture} & \raisebox{0.7cm}{ \begin{tabular}{l}Exits: $CD\subseteq F_k, \ \ \ EF\subseteq F_l,$ \\ $ABCD$ has exits also in $F_i$ and $F_j$. \\ \\$CD \subseteq \mathcal{P}_{x_l=1} \cap F_k.$ \\ \end{tabular} }\\
\hline
\raisebox{-0.1cm}{\begin{tikzpicture}[scale=0.8] \decskeleton \node [above] at (0.5,0.5){3H}; \draw [fill] (0,0) circle [radius=.1];
\draw [ultra thick] (1,0.3)--(1,-0.3);  \end{tikzpicture} } &\begin{tikzpicture}[scale=0.7] \dmH \end{tikzpicture} &\raisebox{0.7cm}{\begin{tabular}{l}Exits: $CE\subseteq F_k, \ \ \ DE\subseteq F_l,$ \\ $ABCD$ has exits also in $F_i$ and $F_j$. \\ \\$CD \subseteq \mathcal{P}_{x_k+x_l=1}.$\\ \end{tabular}}\\
\hline 
\multicolumn{3}{c}{}\\
 \hlineB{3.5}
\multicolumn{3}{!{\vrule width 2pt}c!{\vrule width 2pt}}{Motifs of families of lines}\\
\hlineB{3.5}
Primal Motif & Dual Motif & Necessary conditions \\
\hline
\begin{tikzpicture}[scale=0.8] \decskeleton \node [above] at (0.5,0.5){3I}; \draw [fill] (0,0) circle [radius=.1];\end{tikzpicture} &\begin{tikzpicture}[scale=0.7] \dmI \end{tikzpicture}  &\raisebox{0.9cm}{\begin{tabular}{l}Exits: $CD\subseteq F_k \cap F_l,$ \\ $ABCD$ has  exits also in $F_i$ and $F_j$. \\\end{tabular}}\\
\hline
\begin{tikzpicture}[scale=0.8] \decskeleton \node [above] at (0.5,0.5){3J}; \draw [fill] (-0.3,0.3) circle [radius=.1];
\draw [fill] (-0.3,-0.3) circle [radius=.1];
\draw [ultra thick] (0,0.3)--(-0,-0.3);
\draw [thick, dashed] (-0.3,0.3)--(-0.3,-0.3); \end{tikzpicture} &\begin{tikzpicture}[scale=0.7] \dmJ \end{tikzpicture} &\raisebox{0.7cm}{\begin{tabular}{l}Exits: $BC\subseteq F_i \cap F_j$ \\
$\qquad \ \ \ DE \subseteq F_k \cap F_l$ \\ \\ $AD \subseteq \mathcal{P}_{x_j=1}, AE \subseteq \mathcal{P}_{x_i=1}.$ \\ \end{tabular}}\\
\hline
 \end{tabular}
\end{center}
\caption{The occurrences of motifs for general smooth cubic surfaces.  In the dual motifs the bold lines indicate edges with required exits.}\label{tab:gen_comblines}
\end{table} 



\begin{proposition}\label{prop:inj}
Let $X$ be a smooth tropical surface of degree 3. Suppose $L\sub X$ has motif other than \emph{3I} and \emph{3J}. If $L'\sub X$ is any tropical line, we have
$$\mc_X^\vee(L')=\mc_X^\vee(L)\;\Longrightarrow \; L'=L.$$
\end{proposition}
Alternatively, it  can be formulated as follows: Let $\kr$ be a motif of type either 3A, 3B, 3C, 3D, 3E, 3F, 3G or 3H. Then there is either none or exactly one tropical line on $X$ with dual motif $\kr$. 

\begin{proof}
This is a consequence of \cite[Proposition 8.3]{Vigeland10}, which states that if $\deg X\geq 3$, then any $L\sub X$ not belonging to a two-point family on $X$ is uniquely determined by its {\em set of data}, $\mathcal{D}_X(L)$, introduced in \cite{Vigeland10}. The main difference between $\mc_X^\vee(L)$ and $\mathcal{D}_X(L)$ is that the latter includes the combinatorial type of $L$. However, one can check that if $L\sub X$ has any general motifs (and $\deg X\geq 3)$, then its combinatorial type - and thus $\mathcal{D}_X(L)$ - is uniquely determined by $\mc_X^\vee(L)$ and $X$. To prove the lemma, we argue as follows: Let $L$ be as in the statement, and suppose $\mc_X^\vee(L')=\mc_X^\vee(L)$ for some $L'\sub X$. This implies that $\mathcal{D}_X(L)=\mathcal{D}_X(L')$, and thus, by \cite[Proposition 8.3]{Vigeland10}, that either $L=L'$, or that $L$ belong to a two-point family on $X$. As in the proof of  \cite[Proposition 8.3]{Vigeland10} we check that the latter possibility cannot happen.
\end{proof}

The classification in the last section can be used to determine bounds on the number of  tropical lines on smooth tropical surfaces. More precisely, let $X$ be any smooth tropical surface of degree at least $3$. First, check whether $\ks_X$ contains a motif  3I or 3J. If it does, then by Proposition \ref{prop:inj} $X$ contains infinitely many tropical lines.

Suppose $\ks_X$ contains no motifs 3I or 3J. 
Then Proposition \ref{prop:inj} implies that there is the following injective map: 
\begin{equation*}
\{\text{tropical lines on $X$ of motif $\mathcal{M}$}\}\;\xrightarrow{\phantom{A}\mc_X^\vee\phantom{A}}\;\{\text{occurrences of motif $\mathcal{M}$ in $\ks_X$}\}.
\end{equation*}
Since on a general smooth $X$, every tropical line has general motif, we have:
\begin{proposition}\label{prop:ineq}
Let $\ks$ be a regular unimodular triangulation of $d \, \Delta_3$ without motif of type \emph{3I} and \emph{3J}, where $d\geq 3$. If $X$ is a general smooth tropical surface with subdivision $\ks$, then
\begin{equation*}
\sharp\{\text{tropical lines on $X$}\}\leq \sharp\{\text{general motifs  of $\ks_X$ }\}.
\end{equation*}
\end{proposition}

\begin{remark}
Proposition \ref{prop:ineq} gives a computationally accessible upper bound for the number of tropical lines on a general tropical surface with given subdivision. Namely, if $\ks$ is a subdivision of $d\, \Delta_3$, its general motifs can be found in the following easily programmable way: For each type, identify all occurrences of motifs in $\ks$ by looking for the required cell structures and checking the necessary conditions, as given in Table \ref{tab:classif}. 
\end{remark}

The upper bound given in Proposition \ref{prop:ineq} is not sharp in general. However, in concrete examples, it is possible to improve the inequality, or even find the exact number of tropical lines. We give a detailed example of this in section \ref{sec:cubic}. 

\section{The honeycomb triangulation}\label{sec:cubic}
\begin{figure}[tbp]
\centering
\includegraphics[height=7cm]{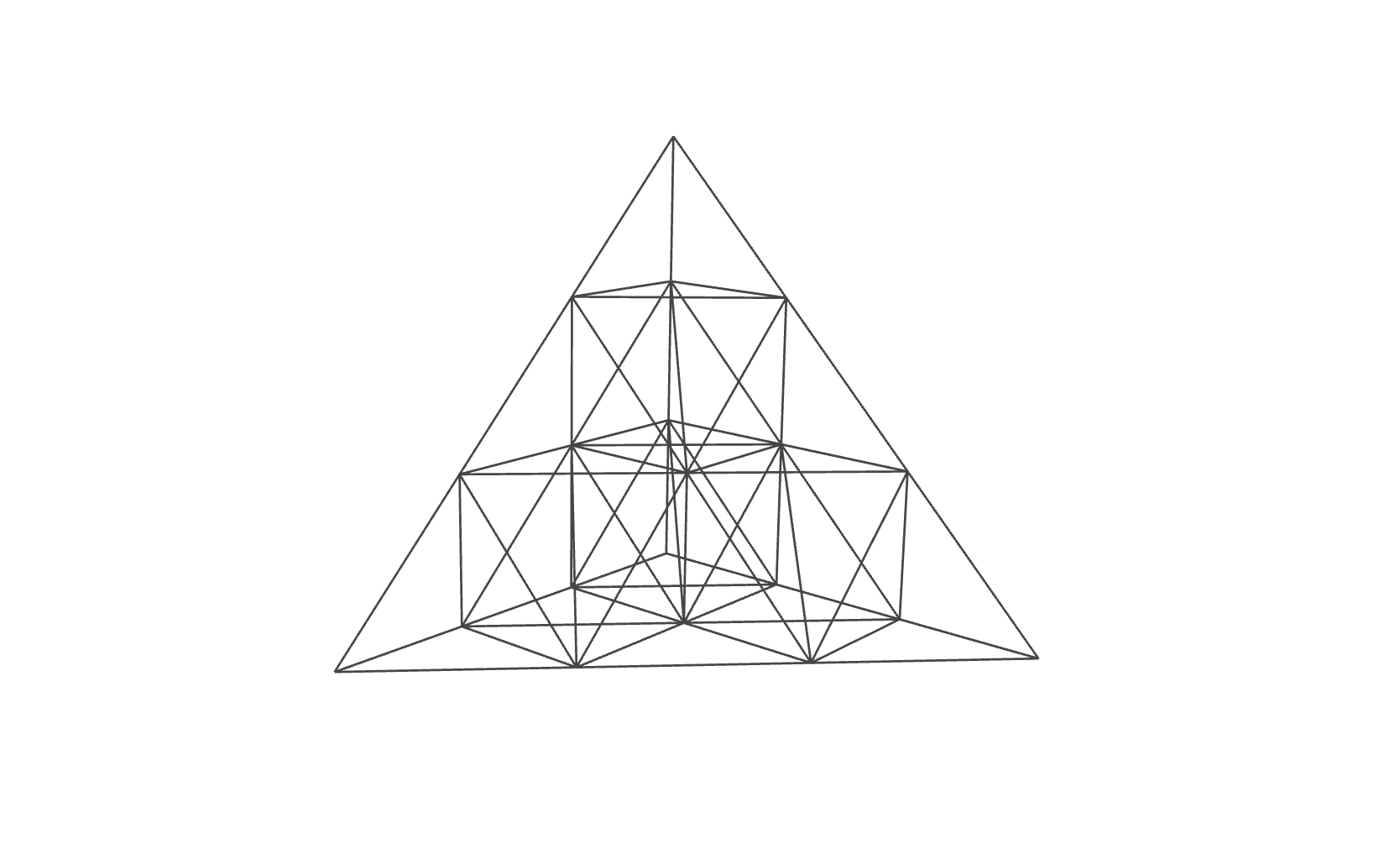}
\caption{The subdivision $\mathcal{H}$, invariant under involutions $(12)$, $(30)$ and $(13)(20)$.}\label{subex}
\end{figure}

%
%
In this section we carefully analyze an example. We fix a dual triangulation of $3 \, \Delta_3$ and we look at the tropical smooth cubic surfaces dual to it  and their lines. We begin by fixing some notation. We label $P_{ijk}$ the lattice point $[i,j,k]$ of $3 \, \Delta_3$. 
We consider the honeycomb triangulation $\mathcal{H}$ shown in Figure \ref{subex} computed with \polymake by the following height function: 
\[ \alpha(x,y,z) =   2 x^2 + 2y^2 + 2z^2 +xy+2xz+2yz. \]

The aim is to prove the following theorem: 

\begin{theorem} \label{theo:ex}
 Let $\mathcal{H}$ be the honeycomb triangulation defined above. 
 \begin{enumerate}
\item A general tropical surface with triangulation $\mathcal{H}$ contains exactly $27$ tropical lines. 
\item Any tropical surface with triangulation $\mathcal{H}$ contains at least $27$ tropical lines. 
\item There exist tropical surfaces with triangulation $\mathcal{H}$ containing infinitely many tropical lines.
\end{enumerate}
\end{theorem}

We will show this through a series of lemmas. The choice of such triangulation is motivated by its symmetries. In fact, we will frequently use that $\mathcal{H}$ is invariant under the subgroup $H\sub S_4$ generated by the three involutions $(12)$, $(30)$ and $(13)(20)$. In particular, $|H|=8$. 

Every surface $X$ we consider in this section  is assumed to have dual subdivision $\mathcal{H}$. Thus $X$ corresponds to a point $c = (c_{000},c_{001},\dotsc,c_{300})$ in the secondary cone $\Sigma(\mathcal{H})$, where the labeling is chosen such that $c_{ijk}$ is the lifting value of $P_{ijk}$. In other words, $X=T(f)$, where 
\[ 
\begin{array}{rl}
f(x,y,z)=&c_{000} + c_{001}z +c_{002}z^2 +  c_{003} z^3 + c_{010} y +  c_{011} yz + c_{012}yz^2  + \\
& +  c_{020} y^2+  c_{021}y^2z + c_{030} y^3 + c_{100} x + c_{101} xz +  c_{102} xz^2  + c_{110} xy  \,+ \\
&+ c_{111} xyz + c_{120} xy^2 + c_{200} x^2 +  c_{201} x^2z+ c_{210} x^2y+ c_{300} x^3. \\
\end{array}
\]
\begin{lemma}\label{lem:ex0}
Let $X$ be a tropical surface with dual triangulation $\mathcal{H}$. Then $X$ has no tropical lines with  motifs \emph{3C},  \emph{3G},  \emph{3H},  \emph{3I} or \emph{3J}.
\end{lemma}
\begin{proof}
It is enough to observe that $\mathcal{H}$ has no occurrences of dual motifs  3C, 3G, 3H,  3I or 3J. This an be checked using the cell structures given in Table \ref{tab:gen_comblines}.
\end{proof}

\begin{lemma}\label{lem:3a3d} Let $X$ be a tropical surface with dual triangulation $\mathcal{H}$. 
\begin{enumerate}
\item[a)] A general $X$ has exactly 12 tropical lines with motif \emph{3A} or \emph{3D}. 
\item[b)] Exactly 4 of the tropical lines in \emph{a)} specialize to a two-point family.
\item[c)] Any $X$ has exactly 12 tropical lines which deforms into motif  \emph{3A} or \emph{3D}. Neither of these deforms into any other motif. 
\end{enumerate}
\end{lemma}
\begin{proof}
\emph{a)} Consider the three occurrences of subcomplexes of dual motifs $\kr_{ \text{3A}},\kr_{\text{3D}},\kr_{\text{3D}}'\sub\mathcal{H}$ shown in Figure~\ref{fig:3a3d}. In $\mathcal{H}$ we find eight subcomplexes  of dual motifs 3A; these are all equivalent modulo $H$ to $\kr_{\text{3A}}$. Furthermore, there are twelve subcomplexes of motifs 3D. Of these, eight are equivalent to $\kr_{\text{3D}}$, while the remaining four are equivalent to $\kr_{\text{3D}}'$.
\begin{figure}[tbp]
\begin{minipage}[b]{.33\linewidth}\begin{center}
\includegraphics[height=4cm]{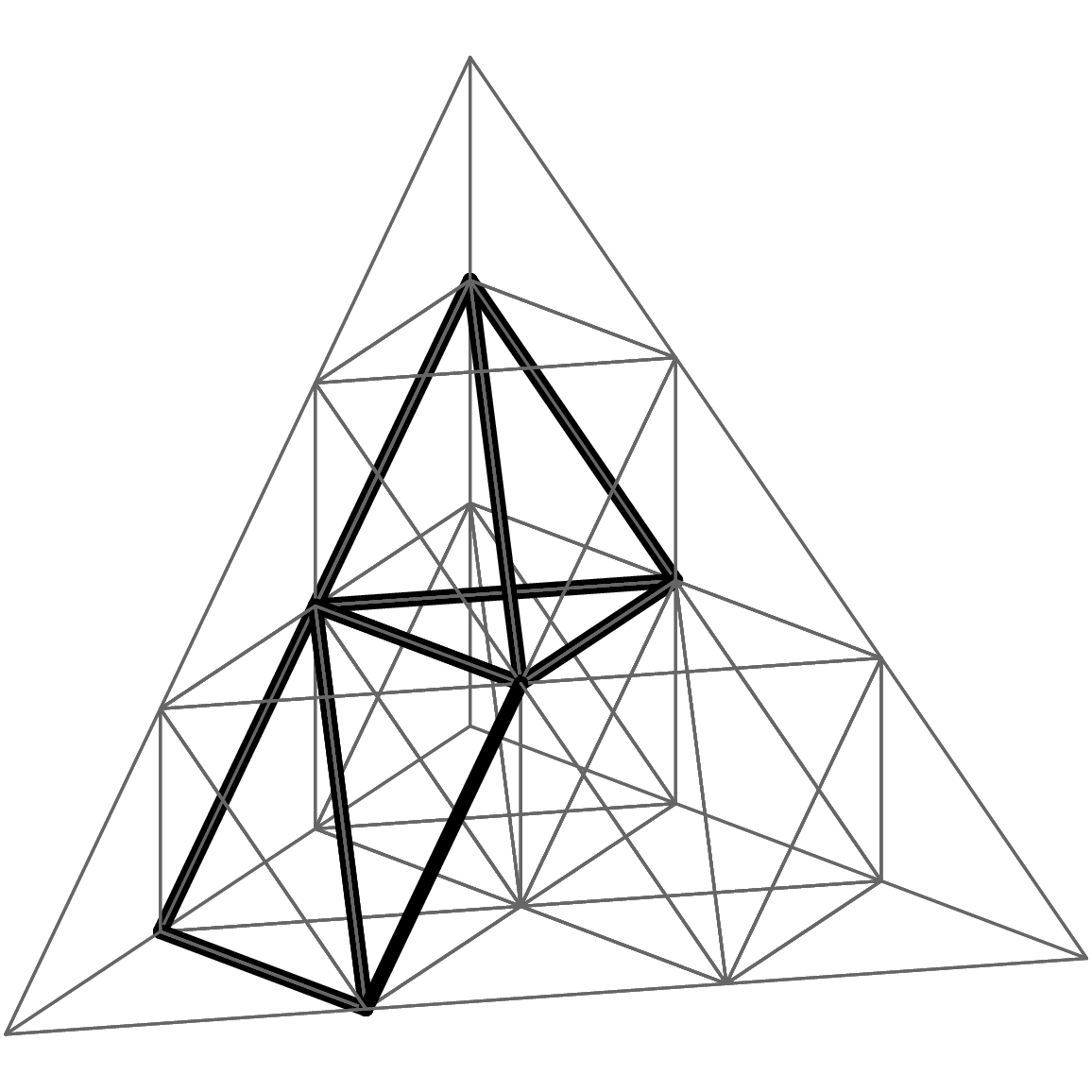}\\
$\kr_{\text{3A}}$
\end{center}
\end{minipage}\hfill
\begin{minipage}[b]{.33\linewidth}\begin{center}
\includegraphics[height=4cm]{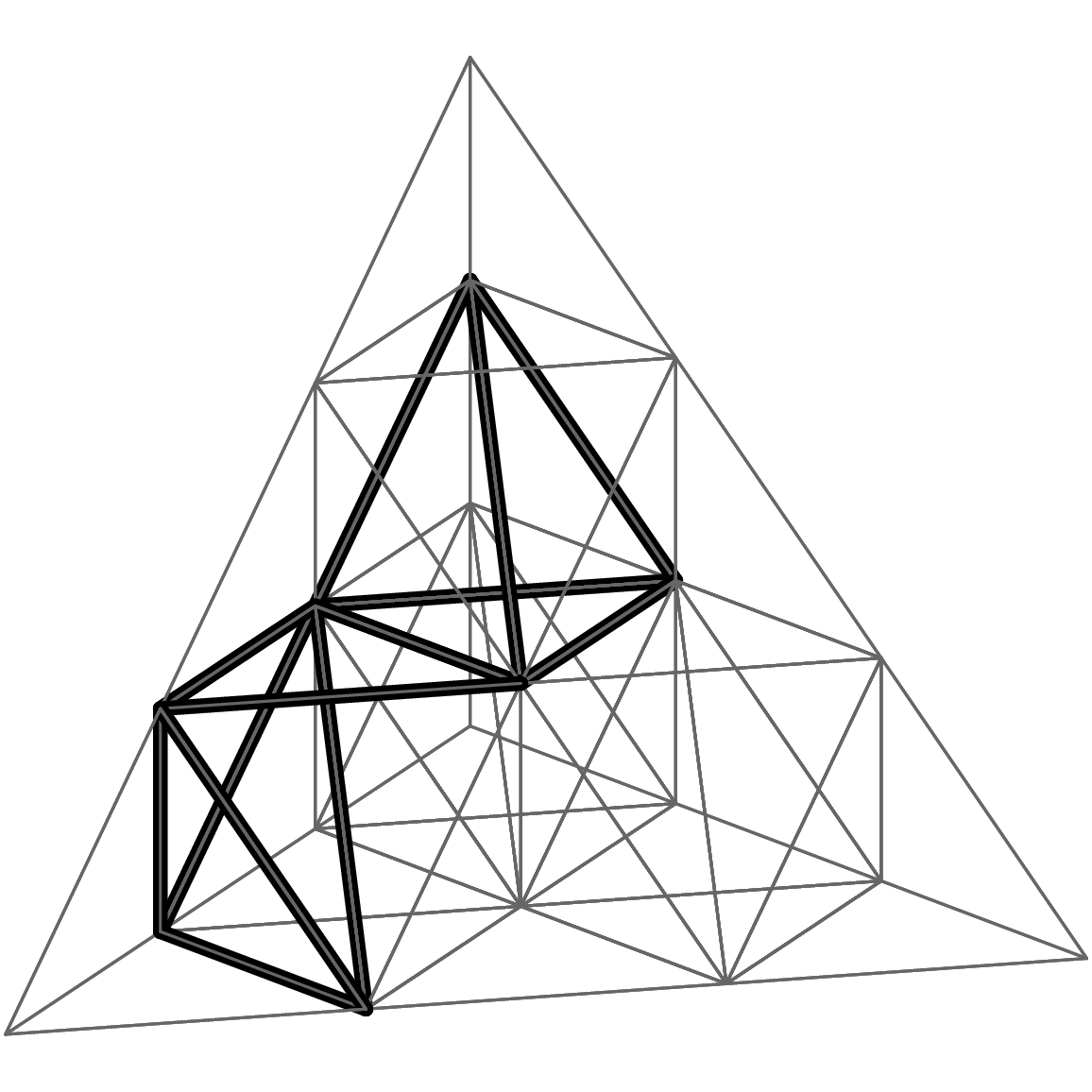}\\
$\kr_{\text{3D}}$
\end{center}
\end{minipage}\hfill
\begin{minipage}[b]{.33\linewidth}\begin{center}
\includegraphics[height=4cm]{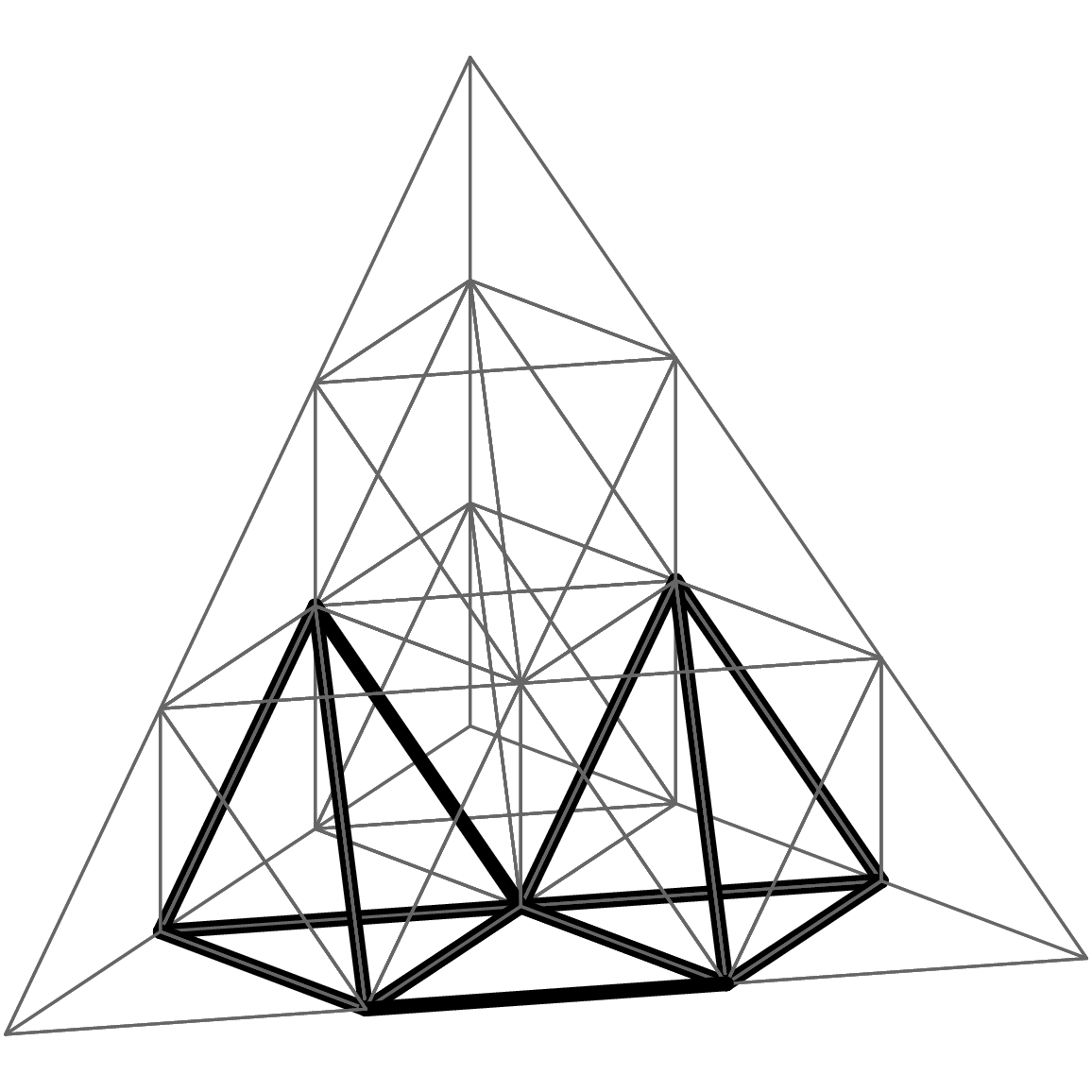}\\
$\kr_{\text{3D}}'$
\end{center}
\end{minipage}\hfill
\caption{Subcomplexes in  $\mathcal{H}$ of dual motifs  3A and 3D}\label{fig:3a3d}
\end{figure}

Let $h_1=c_{210}+c_{002}-c_{201}-c_{011}$, and $h_2=2c_{210}-2c_{120}+c_{020}-c_{200}$. We claim that:
\begin{enumerate}[$i)$]
\item $h_1>0\Longleftrightarrow \text{$\kr_{\text{3A}}$ can be realized uniquely on $X$ as a tropical line with motif 3A}$,
\item $h_1<0\Longleftrightarrow \text{$\kr_{\text{3D}}$ can be realized uniquely on $X$ as a tropical line with motif 3D}$,
\item $h_2\neq0\Longleftrightarrow \text{$\kr_{\text{3D}}'$ can be realized uniquely on $X$ as a tropical line with motif 3D}$.
\end{enumerate}

To prove claims \emph{i)} and \emph{ii)} we sketch the 2-cells of $X$ dual to the three edges $P_{101}P_{111}$, $P_{101}P_{201}$ and $P_{101}P_{210}$ (see Figure \ref{fig:3aor3d}). In the figure, $P$ and $Q$ are the points dual to the tetrahedra $P_{002}P_{011}P_{101}P_{111}$ and $P_{101}P_{111}P_{201}P_{210}$ respectively. 
By Lemma \ref{lem:tresconv}, $X$ contains a line segment trespassing through $P$, with direction vector $\omega_1$. This segment can be extended uniquely to a ray $\ell_1\sub X$, starting somewhere on the polygonal arc $Q'QQ''$. 

A calculation shows that the coordinates of $P$ and $Q$ are
\begin{equation*}
  \begin{split}
    P&=(c_{011}-c_{111}, \ c_{101}-c_{111}, \ c_{101}-c_{111}+c_{011}-c_{002})\\
Q&=(c_{101}-c_{201}, \ c_{101}-c_{111}, \ c_{101}-c_{201}+c_{210}-c_{111}).
  \end{split}
\end{equation*}
In particular, $Q_z-P_z=h_1$. Suppose first that $h_1>0$. Then $Q_z>P_z$, so $\ell_1$ starts in the interior of $QQ''$. Observe that $(QQ'')^\vee= P_{101}P_{111}P_{210}$ has an exit in the direction $\omega_0$, and that $(RR')^\vee= P_{101}P_{200}P_{210}$ has exits in both directions $\omega_2$ and $\omega_3$.  Figure~\ref{fig:3aor3d} shows that if $h_1>0$, then $\ell_1$ can be extended uniquely to a tropical line $L\sub X$ of combinatorial type $((23)(10))$, with one vertex in each of $\inter(QQ'')$ and $\inter(RR')$. The line $L$ has motif 3A, and $\mc_X^\vee(L)=\kr_{\text{3A}}$, so claim \emph{i)} is proved.
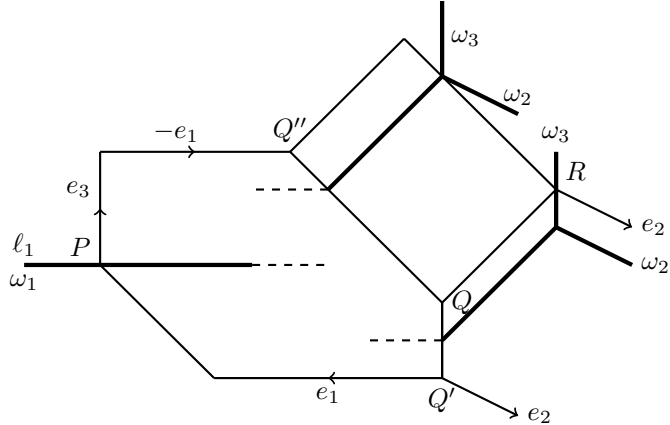
\begin{figure}[ht]
 \centering
 \begin{tikzpicture}
 \begin{scope}[thick,decoration={
    markings,
    mark=at position 0.5 with {\arrow{>}}}
    ] 
    \draw[postaction={decorate}]  (0,0)--(-3,0);
       \draw[postaction={decorate}] (-4.5,1.5) -- (-4.5,3);
    \draw[postaction={decorate}]  (-4.5,3)--(-2,3);
  \end{scope}
      \draw[thick] (0,0)--(0,1);
    \draw[thick]  (-2,3)--(-0.5,4.5);
    \draw[thick]  (-0.5,4.5)--(1.5,2.5);
        \draw[thick]  (1.5,2.5)--(0,1);
            \draw[thick]  (0,0)--(0,1);
    \draw[thick]  (-2,3)--(0,1);
\draw[thick] (-3,0)--(-4.5,1.5);
\draw [ultra thick] (-5.5,1.5)--(-2.5,1.5);
\draw [ultra thick] (0,0.5)--(1.5,2);
\draw [thick, dashed] (-2.5,1.5)--(-1.5,1.5);
\draw [ultra thick] (1.5,2)--(1.5,3);
\draw [ultra thick](1.5,2)--(2.5,1.5);
\draw [thick, dashed](0,0.5)--(-1,0.5);
\draw [ultra thick](-1.5,2.5)--(0,4);
\draw [thick, dashed](-1.5,2.5)--(-2.5,2.5);
\draw [ultra thick](0,5)--(0,4);
\draw [ultra thick] (1,3.5)--(0,4);
\draw [thick,->] (0,0)--(1,-0.5);
\draw [thick,->] (1.5,2.5)--(2.5,2);
\node[above] at  (1,3.5) {$\omega_2$};
\node [right] at (1,-0.5) {$e_2$};
\node [right] at (2.5,1.5) {$\omega_2$};
\node [right] at (2.5,2) {$e_2$};
\node [above] at (1.5,3) {$\omega_3$};
\node [right] at (0,4.5) {$\omega_3$};
\node [above] at (-5.5,1.5) {$\ell_1$};
\node [below] at (-5.5,1.5) {$\omega_1$};
\node [above] at (-2,3) {$Q''$};
\node [below] at (0,0) {$Q'$};
\node [above left] at (-4.5,1.5) {$P$};
\node [right] at (0,1) {$Q$};
\node [above right] at (1.5,2.5) {$R$};
\node [below] at (-1.5,0) {$e_1$};
\node [above] at (-3.5,3) {$-e_1$};
\node [left] at (-4.5,2.5) {$e_3$};
\end{tikzpicture}
    \caption{If the ray $\ell_1$ meets the interior of the segment $QQ''$ (resp. the interior of $QQ'$), it can be extended uniquely to a tropical line on $X$ with motif 3A (resp. 3D).}
  \label{fig:3aor3d}
\end{figure}

Similarly, if $h_1<0$, then $\ell_1$ starts in $\inter(QQ')$. From the facts that $(QQ')^\vee= P_{101}P_{111}P_{201}$ has an exit in the direction $\omega_0$, and that the vertex $R$ allows a trespassing ray with direction $\omega_3$ (cf. Lemma \ref{lem:tresconv}), we see that $\ell_1$ can be extended uniquely to a tropical line on $X$ of combinatorial type $((23)(10))$, with one vertex in $\inter(QQ')$ and the other in the interior of the 2-cell $(P_{101}P_{201})^\vee$. The motif of this line is 3D, and the associated line subcomplex in $\ks_X$ is precisely $\kr_{\text{3D}}$. Thus claim \emph{ii)} is proved.
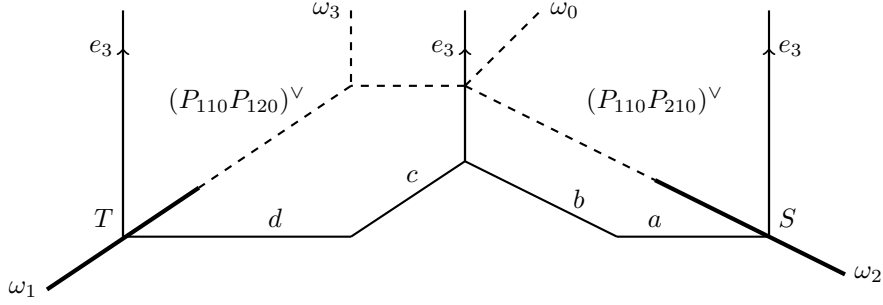
\begin{figure}[h]
\begin{tikzpicture}
\draw [thick] (0,0)--(3,0);
\draw [thick] (3,0)--(4.5,1);
\draw [thick] (4.5,1)--(6.5,0);
\draw [thick] (6.5,0)--(8.5,0);
\draw [thick] (0,0)--(0,3);
\draw [thick] (4.5,1)--(4.5,3);
\draw [thick] (8.5,0)--(8.5,3);
\draw[thick,->] (0,0) -- (0,2.5);
\draw[thick,->] (4.5,1) -- (4.5,2.5);
\draw[thick,->] (8.5,0) -- (8.5,2.5);
\draw [ultra thick] (7,0.75) -- (9.5,-0.5);
\draw [thick, dashed] (4.5,2) -- (7,0.75);
\draw [thick,dashed] (3,2)--(1, 0.65);
\draw [ultra thick] (-1,-0.7) -- (1,0.65);
\draw [thick, dashed] (4.5,2) -- (3,2);
\draw [thick, dashed] (3,3) -- (3,2);
\draw [thick,dashed] (4.5,2) -- (5.5,3);

\node [left] at (0,2.5) {$e_3$};
\node [left] at (4.5,2.5) {$e_3$};
\node [right] at (8.5,2.5) {$e_3$};
\node [above left] at (0,0) {$T$};
\node [above right] at (8.5,0) {$S$};
\node [above] at (2,0) {$d$};
\node [below left] at (4,1) {$c$};
\node [left] at (6.2,0.5) {$b$};
\node [above] at (7,0) {$a$};
\node [left] at (3,3) {$\omega_3$};
\node [right] at (5.5,3) {$\omega_0$};
\node [right] at (9.5,-0.5) {$\omega_2$};
\node [left] at (-1,-0.7) {$\omega_1$};
\node[above] at (1.5,1.5) {$(P_{110}P_{120})^\vee$};
\node[above] at (7,1.5) {$(P_{110}P_{210})^\vee$};
 \end{tikzpicture}
  \caption{A tropical line whose line subcomplex is $\kr_{\text{3D}}'$.}
  \label{fig:3d_real}
 \end{figure}

For claim \emph{iii)} we refer to Figure \ref{fig:3d_real}, showing the 2-cells dual to the edges $P_{110}P_{210}$ and $P_{110}P_{120}$. If the side lengths $a+b\neq c+d$, then $X$ contains a unique tropical line $L$ containing the vertices $S,T\in X$: If $a+b<c+d$, as in Figure \ref{fig:3d_real}, then $L$ has combinatorial type $((13)(20))$ and one vertex in $\inter(P_{110}P_{120})^\vee)$. If $a+b>c+d$ then $L$ has combinatorial type $((10)(23))$ and one vertex in $\inter((P_{110}P_{210})^\vee)$. In both cases, $L$ has one vertex on the edge $(P_{110}P_{210}P_{120})^\vee)$ joining the two 2-cells. The line, $L$ has motif  3D, and $\mc_X^\vee(L)=\kr_{\text{3D}}'$. Furthermore, calculating the vertex coordinates, one finds that $a+b-c-d=h_2$. This proves claim \emph{iii)}. 

Observe that claim \emph{i)} remains valid if we exchange $h_1$ and $\kr_{\text{3A}}$ by $\sigma(h_1)$ and $\sigma(\kr_{\text{3A}})$, where $\sigma$ is any element of $H\sub S_4$, and similarly for the claims  \emph{ii)} and  \emph{iii)}. From this we conclude two things. Firstly, if $c$ lies away from the hyperplanes given by $\sigma(h_1)=0$, for all $\sigma\in H$, then the 16 subcomplexes in the orbits of $\kr_{\text{3A}}$ and $\kr_{\text{3D}}$ give rise to exactly 8 tropical lines on $X_c$. Secondly, if $c$ lies away from the hyperplanes $\sigma(h_2)=0$, for all $\sigma\in H$, then the four subcomplexes in the orbit if $\kr_{\text{3D}}'$ give rise to exactly four tropical lines on $X_c$. Hence, a general $X$ with subdivision $\mathcal{H}$ has exactly $8+4=12$ tropical lines with motif either 3A or 3D.\\[.1cm]
\emph{b)}
Let us first analyze the cases $h_1=0$ and $h_2=0$. If $h_1=0$, then $X$ contains the tropical line with vertices $Q$ and $R$ (see Figure \ref{fig:3aor3d}). It has non-general motif, and it does not belong to any two-point family on $X$.

Next, suppose $h_2=0$. In this case $a+b=c+d$ (cf. Figure \ref{fig:3d_real}), and the lines through $S$ and $T$ with direction vectors $e_2$ and $e_1$ respectively, meet in the point $v:=S-(0,a+b,0)=T-(c+d,0,0)$ on the 1-cell dual to the triangle $P_{210}P_{120}P_{110}$. Since this triangle has exits in both directions $\omega_3$ and $\omega_0$, it follows that $X$ contains the degenerate tropical line with vertex $v$. In fact, it is easy to see that for all $t\geq 0$, the tropical line with vertices $v$ and $v-t(e_1+e_2)$ lies on $X$. Hence $X$ contains the complete two-point family of tropical lines passing through $S$ and $T$.

Now for the specializations. As seen in \emph{a)} the 12 tropical lines in question come in two groups, 8 associated to $\kr_{\text{3A}}$ or $\kr_{\text{3D}}$, and 4 associated to $\kr_{\text{3D}}'$. Suppose $X$ is general, and that $L\sub X$ is in the first group. We can assume that $\mc_X^\vee(L)$ equals either $\kr_{\text{3A}}$ or $\kr_{\text{3D}}$. Any perturbation of $X$ which keeps $h_1\neq 0$, induces a deformation of $L\sub X$ that preserves the motif  of $L$. Hence to obtain a specialization of $L$, we must let $h_1\pil 0$. As observed above, this results in a specialization of $L$ to an isolated tropical line.

Next, let (on a general $X$) $L$ be in the last group, i.e., we can assume $L$ to be the realization of $\kr_{\text{3D}}'$. Choose any perturbation of $X$ such that $h_2\pil 0$. As shown above, this will induce a specialization of $L$ to a degenerate tropical line which belongs to a two-point family on $X$.
\\[.1cm]
\emph{c)} On general $X$, the 12 tropical lines are  found in \emph{a)}. If $X$ is non-general, then either $\sigma(h_1)=0$ or $\sigma(h_2)=0$ for some $\sigma\in H$. Suppose the former. It is enough to consider the case $h_1=0$, in which $X$ contains the tropical line $L_0$ with vertices $Q$ and $R$. As seen in \emph{b)}, $L_0$ is the unique specialization of any realization of $\kr_{\text{3A}}$ or $\kr_{\text{3D}}$. In particular, it deforms into motif 3A and 3D.

We claim that $L_0$ cannot be deformed into any motif other than 3A and 3D. Let $X_0:=X$, and consider any deformation $t\mapsto (L_t,X_t)$ of $L_0\sub X_0$ into some general motif  $\cc$. For each $t$, let $P_t\in X_t$ be the vertex corresponding to $P\in X_0$. Then we know, by Lemma \ref{lem:deftres}, that the $\omega_1$-ray of $L_t$ is trespassing through $P_t$ for each $t$. But this, together with the assumption that the motif  $\cc$ of $L_t$ is general, implies that $\cc$ equals either 3A or 3D. This follows from our discussion in \emph{a)}, in particular Figure \ref{fig:3aor3d}.

Finally, suppose $\sigma(h_2)=0$; as before it is enough to consider the case $h_2=0$. Then $X$ contains the two-point family of tropical lines passing through $S$ and $T$ (cf. Figure \ref{fig:3d_real}). Let $L_{\text{deg}}$ be the degenerate member of this two-point family. It deforms into motif 3D (it is the unique specialization of any realization of $\kr_{\text{3D}}'$), and it is easy to see that it does not deform into any other general motif. As for the non-degenerate tropical lines in the two-point family, none of them have general motif, nor can any of them be deformed into any general motif.

We conclude that the 12 tropical lines found in \emph{a)} all have unique (and distinct) specializations, which satisfy the requirements given in the lemma. This completes the proof.
\end{proof}

\begin{figure}[ht]
\centering
\includegraphics[height=4cm]{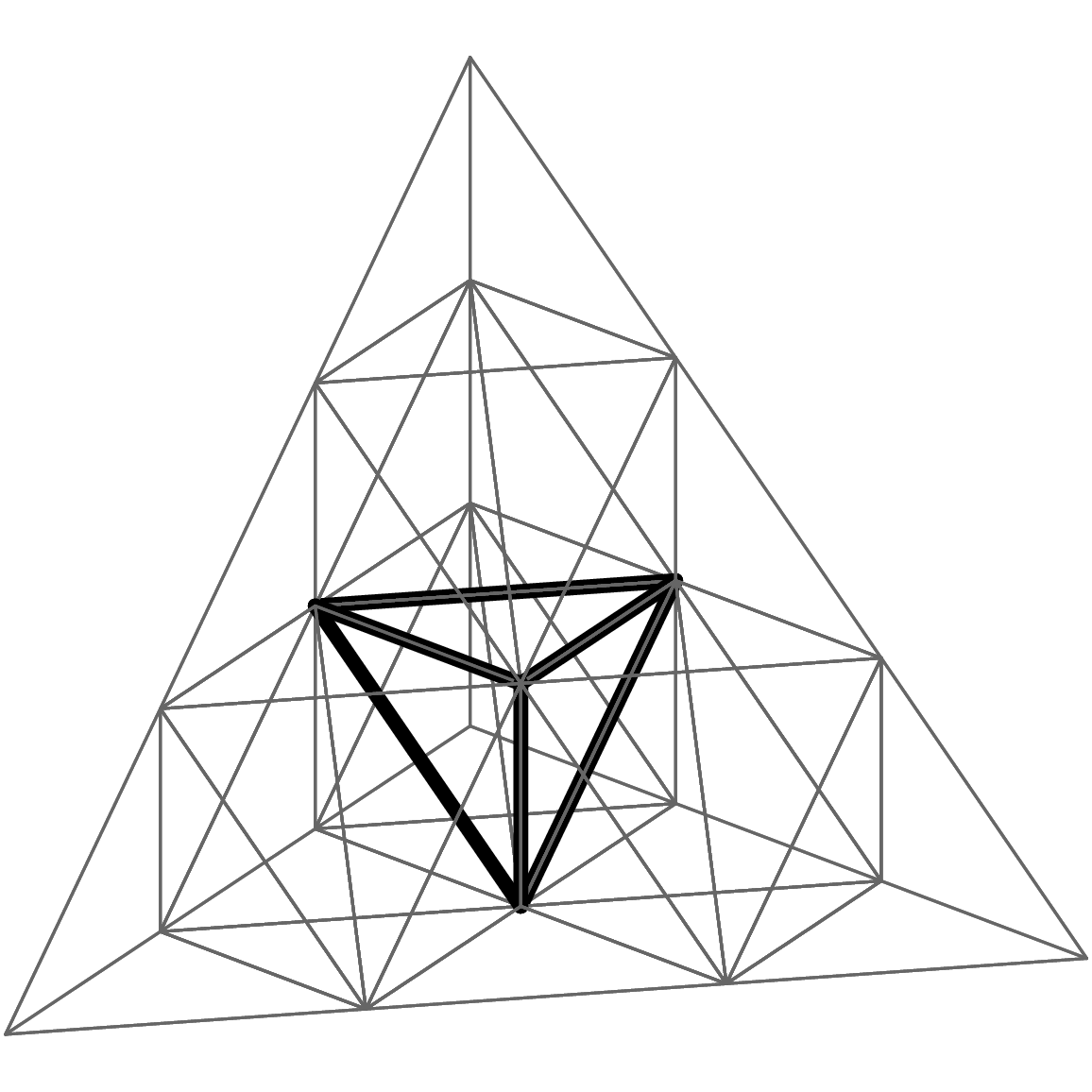}
\caption{The tetrahedron $T$ common to all occurrences of dual motifs 3B.}\label{fig:ex_3b}
\end{figure}

\begin{lemma}\label{lem:ex3b} Let $X$ be a tropical surface with dual subdivision $\mathcal{H}$. 
\begin{enumerate}
\item[a)] A general $X$ has exactly 3 tropical lines with motif \emph{3B}. 
\item[b)] Each of the tropical lines in \emph{a)} specializes to a two-point family.
\item[c)] Any $X$ has exactly 3 tropical lines which deforms into motif \emph{3B}. Neither of these deforms into any other general motif.
\end{enumerate}
\end{lemma}

\begin{proof}
\emph{a)} There are 12 occurrences of  subcomplexe of dual motifs 3B in $\mathcal{H}$. All of these contain the tetrahedron $T$, shown in Figure \ref{fig:ex_3b}, with vertices $P_{011},P_{101}, P_{110}$ and $P_{111}$. Using Lemma \ref{lem:tresconv} we see that the dual vertex $T^\vee\in X$ allows trespassing line segments in three directions simultaneously: $e_1+e_2$, $e_1+e_3$ and $e_2+e_3$. Drawing the shapes of the 2-cells adjacent to $T^\vee$, one sees immediately that each of these trespassing line segments can be extended to a tropical line on $X$. For general $X$, each extension is unique on $X$, and the three resulting tropical lines all have motif $3B$.\\[.1cm]
\emph{b)} Non-generality in this case means that at least one of the three trespassing line segments in \emph{a)} meets a second vertex of $X$. One can check that this always allows for a second trespassing, resulting in a two-point family on $X$. It is even possible for the line segment to meet a third vertex of $X$, giving rise to a 2-dimensional two-point family on $X$. Any of the tropical lines in \emph{a)} specializes to both a 1-dimensional and 2-dimensional family obtained in this way.\\[.1cm]
\emph{c)} For any of the two-point families described in \emph{b)}, none of its members has general motif. Using arguments similar to those in the proof of Lemma \ref{lem:3a3d}, it is not hard to show that there is exactly one tropical line in the family that can be deformed into some general motif, which must be 3B. The truth of the statement follows from this.
\end{proof}

\begin{lemma} \label{lem:ex3e}
Let $X$ be a tropical surface with dual subdivision $\mathcal{H}$. 
\begin{enumerate}  
\item[a)] A general $X$ has exactly 4 tropical lines with motif \emph{3E}. 
\item[b)] Each of the tropical lines in \emph{b)} specializes to a two-point family.
\item[c)] Any $X$ has exactly 4 tropical lines which deforms into  motif \emph{3E}. Neither of these deforms into any other general motif.
\end{enumerate}
\end{lemma}

\begin{figure}[ht]
\centering
\begin{tikzpicture}
\draw [thick] (-1,3)--(2,3);
\draw [thick] (2,3)--(3,2);
\draw [thick] (3,2) --(3,0);
\draw [thick] (1,0) --(3,0);
\draw [thick] (1,0) --(-1,2);
\draw [thick] (-1,3) --(-1,2);
\draw[thick,-<] (-1,3) -- (1,3);
\draw[thick,-<] (3,2) -- (3,0.5);
\draw [ultra thick] (2,2.7) -- (2,4);
\draw [thick, dashed](2,2)--(2,2.7);
\draw [thick, dashed] (0,2)--(2,2);
\draw [ultra thick] (-2,2) -- (0,2);
\draw [thick,dashed] (2,2) -- (3,1);
\draw [thick,dashed] (3,1) -- (4,1.5);
\draw [thick,dashed] (3,1) -- (4,0.5);
\node [above] at (1,1) {$(P_{101}P_{111})^\vee$};
\node [above right] at (2,3) {$P$};
\node [above left] at (2,2) {$v_1$};
\node [above right] at (2,3) {$P$};
\node [ right] at (2,4) {$\omega_3$};
\node [above] at (-2,2) {$\omega_1$};
\node [below] at (-1,2) {$Q$};
\node [below] at (1,0) {$R''$};
\node [below] at (3,0) {$R$};
\node [right] at (3,2) {$R'$};
\node [right] at (3,0.5) {$e_3$};
\node [above] at (1,3) {$e_1$};
\node [left] at (-1,2.5) {$a$};
\node [left] at (0,1) {$b$};
\node [right]  at (2.5,2.5) {$c$};
 \end{tikzpicture}
\caption{A tropical line on $X$ with motif 3E.}
\label{3e_real}
\end{figure}
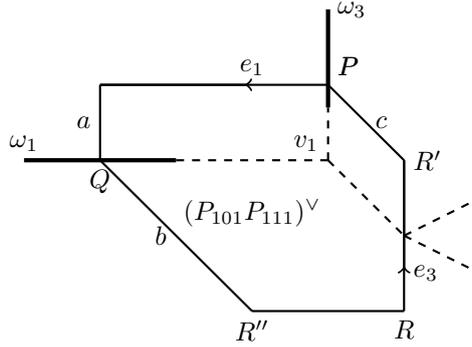

\begin{proof}
There are 8 occurrences of subcomplexes of motif 3E in $\mathcal{H}$, all equivalent modulo $H$. These 8 can be divided into 4 pairs, such that the subcomplexes in each pair contain the same tetrahedra. One of these pairs, $\kr_{\text{3E}}$ and $\kr_{\text{3E}}'$, is shown in Figure \ref{fig:3e}. 

\begin{figure}[ht]
\centering
\includegraphics[height=4cm]{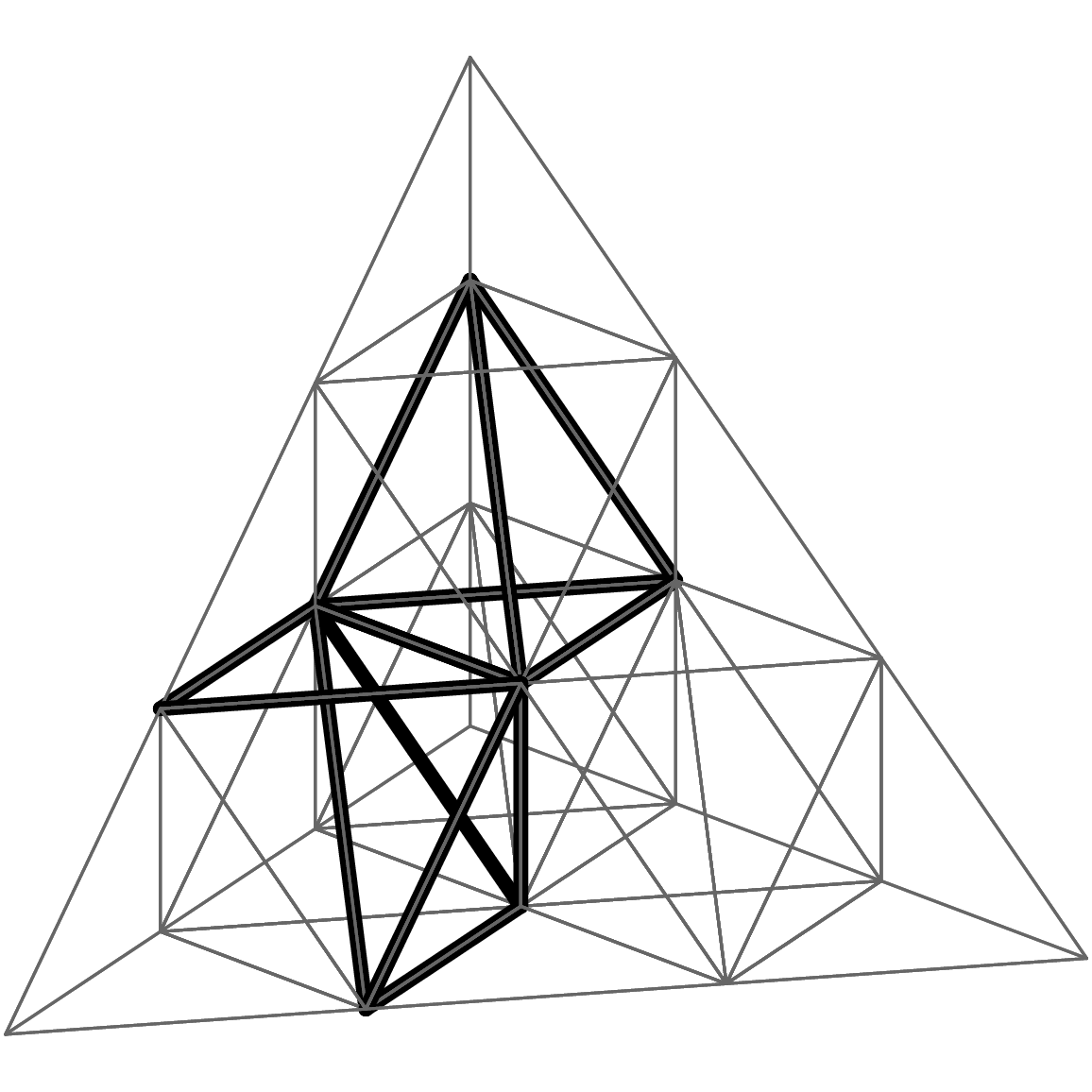}\qquad\qquad\includegraphics[height=4cm]{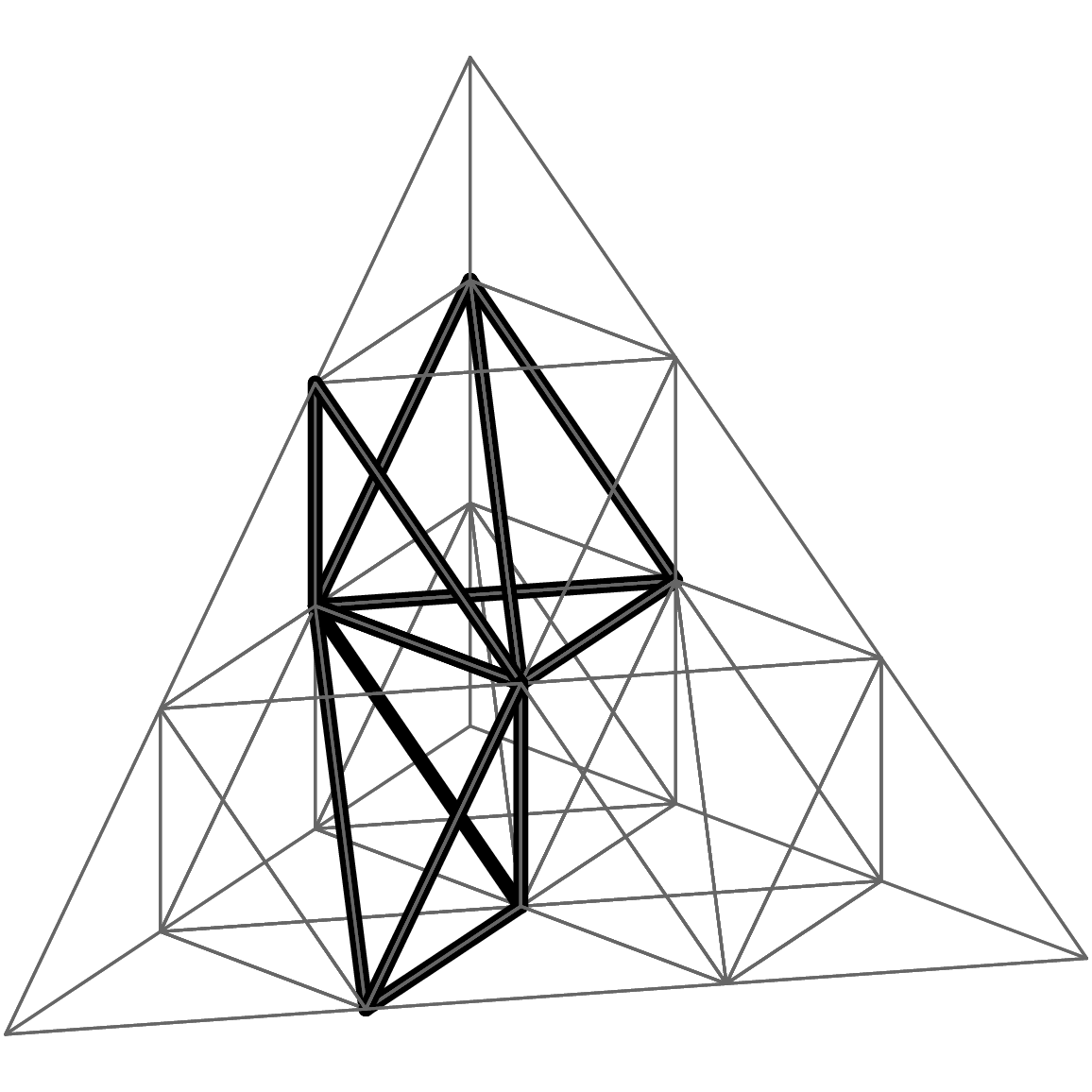}
\caption{Subcomplexes $\kr_{\text{3E}}$ (to the left) and $\kr_{\text{3E}}'$ (to the right).}\label{fig:3e}
\end{figure}

We claim that a general $X$ contains exactly one tropical line $L$ with dual motif  either $\kr_{\text{3E}}$ or $\kr_{\text{3E}}'$. To prove this, we refer to Figure \ref{3e_real}, which shows the 2-cell dual to $P_{101}P_{111}$. The cell is a parallel hexagon whose edge directions are given in the figure. The vertices $P$ and $Q$ are the duals of the tetrahedra $P_{101}P_{110}P_{111}P_{210}$ and $P_{011}P_{101}P_{111}P_{002}$, allowing (by Lemma \ref{lem:tresconv}) trespassing in directions $\omega_3$ and $\omega_1$ respectively. Let $L$ be the tropical line with vertices $v_1=P-(0,0,a)$ and $v_2=v_1-(\min(b,c),0,\min(b,c))$. Observe that $v_2$ lies either on the edge $RR'$ (if $c\leq b$) or on the edge $RR''$ (if $b\leq c$). Hence, since both $(RR')^\vee=P_{111}P_{101}P_{201}$ and $(RR'')^\vee=P_{111}P_{101}P_{102}$ has exits in directions $\omega_2$ and $\omega_0$, this ensures that $L\sub X$. For a general $X$, we can assume $b\neq c$. If $c<b$ (as shown in Figure \ref{3e_real}), we have $v_2\in \inter(RR')$, giving $\mc_X^\vee(L)=\kr_{3E}$. If $b<c$, then $v_2\in \inter(RR'')$, and $\mc_X^\vee(L)=\kr_{\text{3E}}'$. In either case it is clear that $L$ is the only tropical line on $X$ passing through $P$ and $Q$. This proves the claim.

By symmetry, the same argument applies to the three other pairs of dual motifs 3E, giving a total of 4 lines with motif  3E.

Parts \emph{b)} and \emph{c)} are proved in a similar fashion as in the corresponding parts of Lemma \ref{lem:3a3d}.
\end{proof}

\begin{lemma}\label{lem:ex3f}
Any tropical surface $X$ with subdivision $\mathcal{H}$ has exactly 8 tropical lines with motif 3F. Neither of these specializes into any other motif.
\end{lemma}

\begin{figure}[h]
  \centering
\begin{tikzpicture}
 \begin{scope}[thick,decoration={
    markings,
    mark=at position 0.5 with {\arrow{>}}}
    ] 
          \draw[postaction={decorate}]  (4.5,1)--(6.5,0);
          \draw [postaction ={decorate}] (6.5,0)--(8.5,0);
          \draw [postaction= {decorate}] (0,0)--(3,0);
    \end{scope}
    \begin{scope}[thick,decoration={
    markings,
    mark=at position 0.5 with {\arrow{<}}}
    ] 
    \draw[postaction={decorate}]  (3,0)--(4.5,1);
    \end{scope}
\draw [thick] (0,0)--(3,0);
\draw [thick] (6.5,0)--(8.5,0);
\draw [thick] (0,0)--(-1,3);
\draw [thick] (4.5,1)--(5.5,3);
\draw [thick] (4.5,1)--(3.5,3);
\draw [thick] (8.5,0)--(9.5,3);
\draw[thick,->] (4.5,1) -- (5.25,2.5);
\draw[thick,->] (4.5,1) -- (3.75,2.5);
\draw[thick,->] (8.5,0) -- (9.25,2.25);
\draw[thick,->] (0,0) -- (-0.75,2.25);
\draw [thick,dashed] (6.5,1) -- (8,0.25);
\draw [ultra thick] (8,0.25)--(9.5,-0.5);
\draw [ultra thick] (-1,-0.7) -- (0.25,0.15);
\draw [thick, dashed] (0.25,0.15)--(1.5,1);
\draw [ultra thick] (3,1) -- (5.5,1);
\draw [thick, dashed] (1.5,1)--(3,1);
\draw [thick, dashed] (5.5,1)--(6.5,1);
\draw [dashed, thick] (1.5,1) -- (1,2);
\draw [thick, dashed] (6.5,1) -- (7,2);
\node [below] at  (5.5,0.5) {$e_2$};
\node [below] at  (4,0.5) {$e_1$};
\node [left] at (3.75,2.5) {$\omega_0$};
\node [left] at (5.25,2.5) {$e_3$};
\node [right] at (9.25,2.25) {$e_3$};
-\node [left] at (-0.75,2.25) {$\omega_0$};
\node [above left] at (0,0) {$R$};
\node [above right] at (8.5,0) {$P$};
\node [above] at (2,0) {$d$};
\node [below left] at (4,1) {$c$};
\node [left] at (6.2,0.5) {$b$};
\node [above] at (7,0) {$a$};
\node [above]  at (1,2.5) {$(P_{111}P_{120})^\vee$};
\node [above] at (7,2.5) {$(P_{110}P_{210})^\vee$};
\node [right] at  (7,2) {$\omega_3$};
\node [left] at (1,2) {$\omega_0$};
\node [right] at (9.5,-0.5) {$\omega_2$};
\node [left] at (-1,-0.7) {$\omega_1$};
\node [below] at (7.5,0) {$e_2+e_3$};
\node [below] at (1.5,0) {$e_2+e_3$};
 \end{tikzpicture}
 \caption{A tropical line on $X$ with motif $3F$.}\label{fig:3f_real}
 \end{figure}
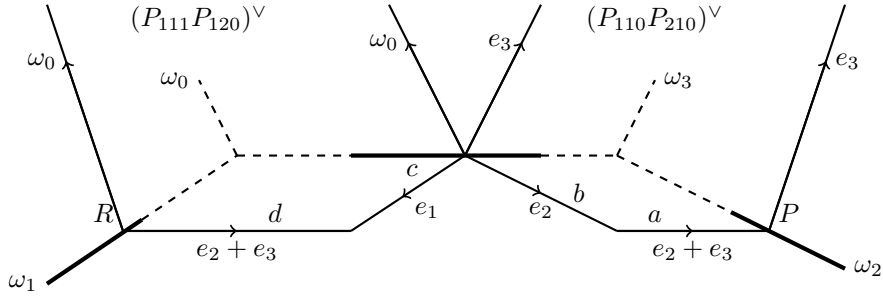

\begin{proof}
Modulo $H$, the only occurrences of subcomplexes of motif 3F in $\mathcal{H}$ are $\kr_{\text{3F}}$ and $\kr_{\text{3F}}'$, shown in Figure \ref{fig:sub3f}. Both have orbits of length 4 under the action of $H$.
\begin{figure}[tbp]
\centering
\includegraphics[height=4cm]{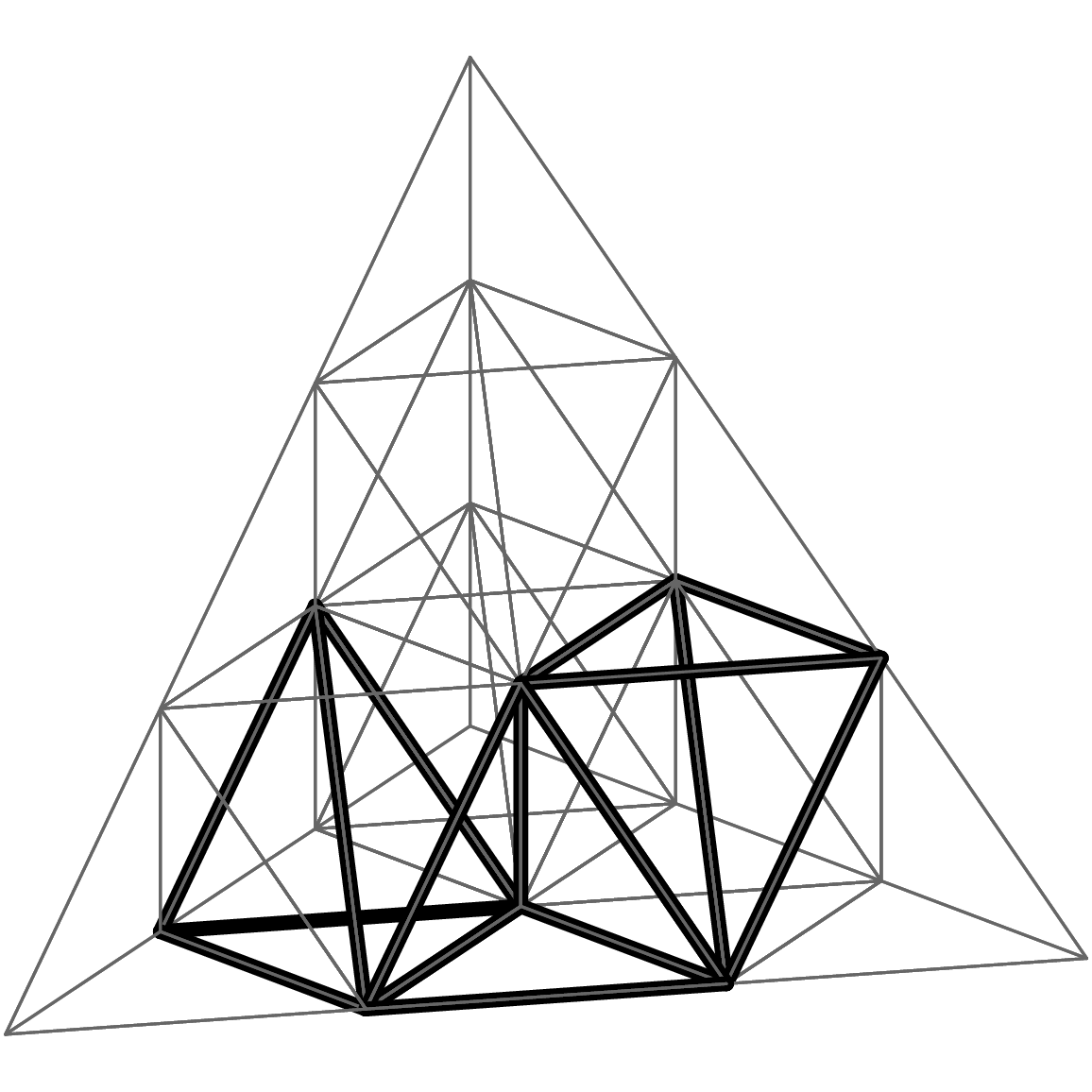}\qquad\qquad\includegraphics[height=4cm]{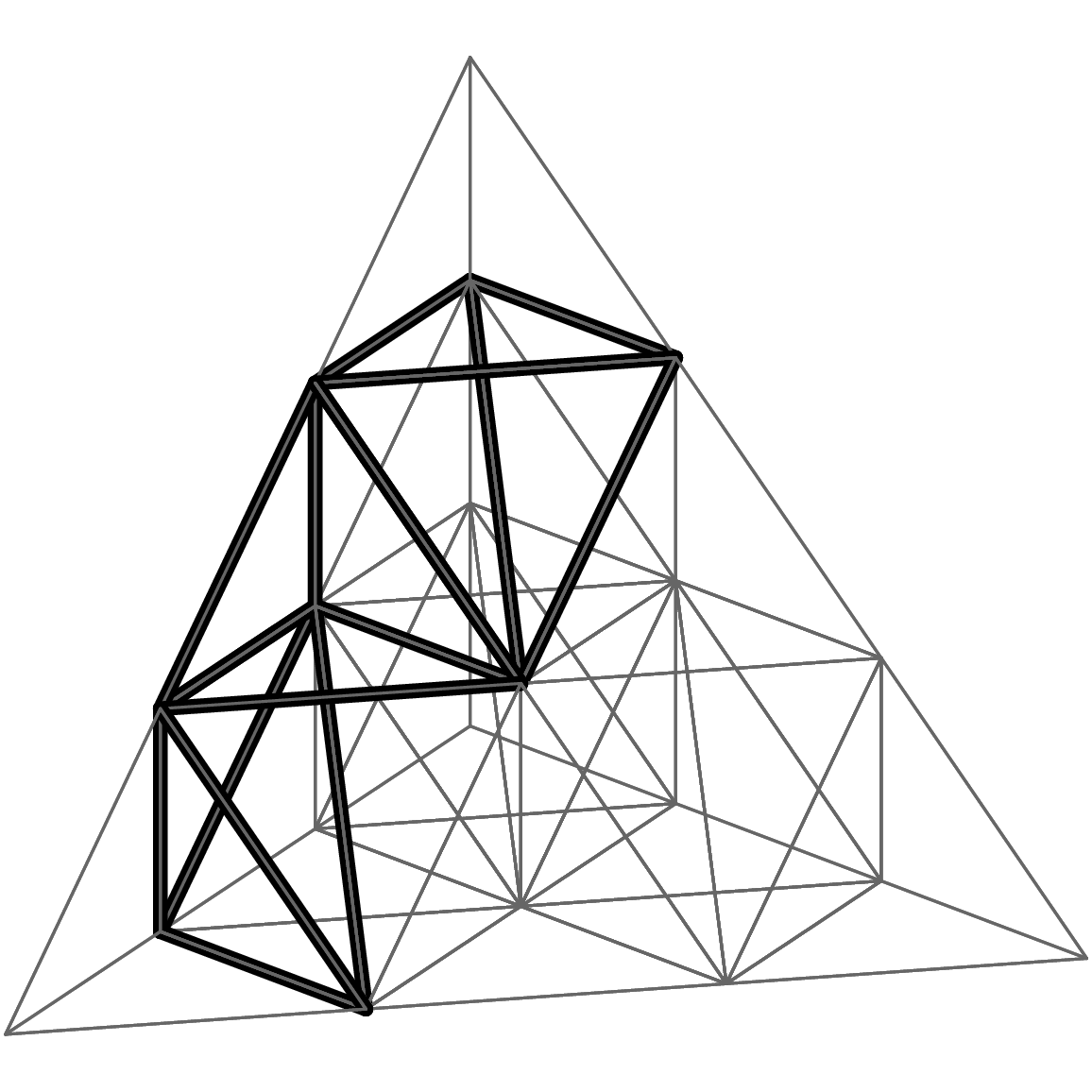}
\caption{Subcomplexes $\kr_{\text{3F}}$ (to the left) and $\kr_{\text{3F}}'$ (to the right).}\label{fig:sub3f}
\end{figure}

It is not hard to see that any $X$ contains exactly one tropical line with line subcomplex $\kr_{\text{3F}}$. Indeed, Figure \ref{fig:3f_real} shows how to construct such a tropical line. For uniqueness we can e.g. apply Lemma \ref{lem:uniqueline}: Denoting the side lengths by $a,b,c,d$, as indicated, we find that $\overrightarrow{RP}=(0,d,d)+(-c,0,0)+(0,b,0)+(0,a,a)=(-c,a+b+d,a+d)$. Since $a,b,c,d$ are strictly positive, the lemma implies that there is a unique tropical line through $P$ and $R$, and, {\em a fortiori}, that there is a unique line on $X$ with associated dual motif $\kr_{\text{3F}}$. The same argument applies to the subcomplexes in the orbit of $\kr_{\text{3F}}$. 

Similarly, by studying the 2-cells dual to $P_{111}P_{102}$ and $P_{101}P_{201}$, one can show that $X$ always contains exactly one tropical line with $\kr_{\text{3F}}'$ as its dual motif. Hence we have a total of 8 tropical lines with motif 3F.
\end{proof}

The lemmas \ref{lem:ex0} through \ref{lem:ex3f} provide everything needed to prove Theorem \ref{theo:ex}.

\begin{proof}[Proof of Theorem \ref{theo:ex}]
\emph{a)} To sum up, we have on a general $X$, 12 tropical lines with motifs 3A or 3D, 3 lines with 3B, 4 lines with 3E, 8 lines with 3F and none with 3C, 3G, 3H, 3I or 3J. No tropical line can have more than one motif on $X$, hence the total number of lines is exactly $12+3+4+8=27$.\\[.1cm]
\emph{b)} Part \emph{c)} of the lemmas \ref{lem:3a3d} through \ref{lem:ex3e}, and Lemma \ref{lem:ex3f} identifies, on any $X$, four sets of tropical lines. Moreover, it follows from the same results that these four sets are mutually disjoint, and contains altogether $27$ tropical lines.\\[.1cm]
\emph{c)} As shown in part  \emph{b)} of the lemmas \ref{lem:3a3d} through \ref{lem:ex3e} there exist tropical surfaces $X$ with subdivision $\mathcal{H}$ containing one or more two-point families of tropical lines. In particular, such $X$ has infinitely many tropical lines.
\end{proof}

\section{Final remarks}

The analysis of tropical surfaces dual to the honeycomb triangulation in Section \ref{sec:cubic} can be carried out for any of the $14\,373\,645$ combinatorial types of smooth cubic surfaces. Implementing such analysis  comes with  interesting algorithmical and computational challenges. First,  it is necessary to enumerate all the occurrences of dual motifs by looking for the right subcomplexes and  their exits. Second, the conditions of realizability of the motifs by lines on a surface needs to be understood. As we have seen in the example, these conditions are different in the various motifs.

Moreover, in this article we have treated each line on the surface individually. The classical geometry of lines on general smooth cubic surfaces is rich of results about the  incidences of the lines. Of course, it is interesting to  look into these results in the tropical setting too. 

We conclude with a question which brings us back to the correspondence between Cayley and Salmon. 
\begin{question}
Does there exist a regular unimodular triangulation $\mathcal{S}$ of $3\, \Delta_3$ such that for every $c$ in the secondary cone $\Sigma(\ks)$  the  tropical surface $X_c$ contains exactly $27$ tropical lines?
\end{question}

\bibliographystyle{plain}

\end{document}